\documentclass{article}
\usepackage{arxiv}
\usepackage[utf8]{inputenc} 
\usepackage[T1]{fontenc}    
\usepackage{hyperref}       
\usepackage{url}            
\usepackage{booktabs}       
\usepackage{amsfonts}       
\usepackage{nicefrac}       
\usepackage{microtype}      
\usepackage{lipsum}
\usepackage{amsmath,amsthm,amsfonts,amssymb,amsbsy,bm,mathtools,textcomp}
\usepackage{graphicx}
\usepackage[rflt]{floatflt}
\usepackage{hyperref}
\usepackage{xcolor}
\usepackage{subcaption}
\usepackage{adjustbox}
\usepackage{tabularx}
\usepackage{array}
\usepackage{siunitx}
\usepackage{comment}
\usepackage[numbers]{natbib}
\usepackage[textsize=tiny]{todonotes}
\usepackage{makecell}
\usepackage[listings]{tcolorbox}
\usepackage{url}
\usepackage{breqn}
\usepackage[normalem]{ulem}
\usepackage{wrapfig}
\usepackage{newunicodechar}
\newunicodechar{，}{,}
\usepackage{comment}
\usepackage[final]{pdfpages}  
\usepackage{xcolor}
\usepackage{soul}
\usepackage{tabularx}
\usepackage{enumitem}
\usepackage{multirow}
\usepackage{lipsum}
\usepackage{algorithm}
\usepackage{algpseudocode}
\usepackage{float}

\title{
Linearized PINN with pretrained nonlinear layers
}

\author{
  Wenhao Chen\\
  Civil and Environmental Engineering Department\\
  University of Illinois Urbana-Champaign\\
  Urbana, IL 61801, USA\\
  \texttt{wenhaoc3@illinois.edu}
  \And  
  Alexandre Tartakovsky\\
  Civil and Environmental Engineering Department\\
  University of Illinois Urbana-Champaign\\
  Urbana, IL 61801, USA\\
  Pacific Northwest National Laboratory\\
  Richland, WA 99352, USA\\
  \texttt{amt1998@illinois.edu}
}

\begin{document}
\maketitle
\begin{abstract}

We propose a Linearized Physics-Informed Neural Network (lPINN), a reduced-order neural basis method for solving forward and inverse partial differential equations. During an offline stage, lPINN learns operator-compatible, continuous neural basis functions from an ensemble of numerical solutions. These functions are pretrained using a dataset consisting of PDE solutions and their derivatives. The learned basis functions are differentiable through automatic differentiation, enabling physics-based inference. For each new problem instance, the basis functions are frozen, and the solution is determined by minimizing the governing-equation residual, together with applicable initial, boundary, regularization, and observational terms. This formulation differs from direct surrogate and operator-learning approaches because the training dataset determines the basis functions in the offline step, whereas the instance-specific solution is obtained by enforcing the governing physics online. 
For linear differential operators, the online lPINN problem reduces to regularized linear least squares. For nonlinear operators, it remains nonlinear but is restricted to the low-dimensional basis coefficients. In inverse problems, the reduced coefficients and unknown physical parameters are estimated jointly from physics residuals and sparse observations. 
Compared to vanilla PINN, pretraining basis functions amounts to estimating parameters in (nonlinear) hidden layers in an offline step and learning parameters in the last (linear) layer in the online step.  

We evaluate lPINN on forward and inverse problems for the advection--diffusion equation, Burgers' equation, and the nonlinear pendulum equation. The experiments compare derivative-matching and residual-based basis pretraining and demonstrate compatibility with Fourier feature networks for oscillatory solutions. Compared with vanilla PINNs, lPINN achieves lower solution and parameter errors in the reported test cases while reducing online inference times from about one order of magnitude to more than three orders of magnitude, with the largest gains generally observed when residual or measurement data are limited. Cross-resolution tests further show that the learned continuous representation can be evaluated on finer meshes without retraining and with nearly unchanged accuracy. These results demonstrate that lPINN can provide an effective alternative to PINN when the offline cost can be amortized over many forward or inverse queries.

\end{abstract}

\section{Introduction}

Physics-informed neural networks (PINNs) provide a flexible framework for forward and inverse problems governed by differential equations by representing the unknown state with a coordinate-based neural network and incorporating governing equations, initial and boundary conditions, and available observations into the training objective~\cite{raissi2019physics,karniadakis2021physics}. In forward problems, the network parameters are optimized to obtain a continuous approximation of the PDE solution that satisfies the prescribed physical constraints. In inverse problems, unknown physical parameters are included among the optimization variables and estimated jointly with the state from physics constraints and observational data~\cite{tartakovsky2020physics,tipireddy2019comparative}. This formulation has enabled applications in fluid and solid mechanics, transport, parameter identification, and constitutive-model discovery. However, its flexibility comes at a substantial computational cost. Conventional PINN training estimates all network parameters for each problem instance, producing a high-dimensional, nonconvex, and often poorly conditioned optimization problem. Convergence can depend strongly on initialization, loss weighting, and sampling, and the network may fail to resolve physically important features. These difficulties are particularly pronounced for stiff, advection-dominated, oscillatory, and multiscale systems~\cite{krishnapriyan2021characterizing}. The PINN computational burden becomes more consequential in problems requiring multiple PDE solves, which limits PINN's applicability in uncertainty quantification, sensitivity analysis, optimization, and time-critical prediction.

A variety of techniques have been introduced to improve the optimization and representation properties of PINNs. Quasi-Newton and related methods can increase the stability and accuracy of training for nonlinear and stiff problems~\cite{urban2025unveiling}. Adaptive weighting and preconditioning strategies seek to balance physical, boundary, initial, and observational terms in the objective~\cite{wang2021understanding}. Fourier features and spectral architectures improve the representation of oscillatory or high-frequency components and can mitigate the spectral bias of standard multilayer perceptrons~\cite{tancik2020fourier,wang2022and}. Although these developments can improve single-instance training, they do not remove the need to optimize a large nonlinear representation for every new problem. The present work addresses this repeated-optimization bottleneck by combining ideas from reduced-order modeling, differentiable neural fields, and physics-informed coefficient inference.

We propose lPINN as a neural reduced-order method for problems requiring solving multiple forward and inverse PDE problems. During the offline stage, an ensemble of numerical solutions is decomposed into a mean field and fluctuations. Coordinate-based neural networks approximate the mean and a shared collection of basis functions for reduced-order representation of the fluctuations. The basis functions may be pretrained by matching states and operator-relevant derivatives or by augmenting state reconstruction with governing-equation, initial-condition, and boundary-condition residuals. During the online stage, the coefficients are determined by minimizing the governing-equation residual together with applicable initial, boundary, regularization, and observational terms. Because the basis functions are differentiable, this inference is performed within the standard PINN framework: automatic differentiation supplies the derivatives of each basis function, while optimization is confined to the coefficients of the last layer. For linear differential operators and linear constraints, the online problem is regularized linear least squares under the usual rank or positive-regularization conditions. For nonlinear operators, it is a reduced-dimensional nonlinear least-squares problem. In inverse problems, the reduced coefficients and unknown physical parameters are estimated jointly.

Standard reduced-order methods are built on the offline--online paradigm for parametrized differential equations. In an offline stage, high-fidelity snapshots are used to construct a low-dimensional space of basis functions, commonly through proper orthogonal decomposition; in the online stage, the state is approximated in this space, and its reduced coordinates are determined from a projected governing model~\cite{hesthaven2016certified,quarteroni2016reduced}. This decomposition can yield substantial savings when the offline cost is amortized over many PDE solves for different parameter queries. Standard bases are nevertheless tied to the spatial or space--time discretization on which the snapshot vectors and reduced operators are defined. Derivatives are inherited from the full-order discretization, and transferring the basis to a different mesh or evaluating the state at arbitrary coordinates generally requires interpolation, reconstruction, or a new discretization. The lPINN method retains the offline--online structure but replaces a discrete snapshot basis with a coordinate-dependent neural basis. The resulting basis functions are continuous and differentiable with respect to space and time through automatic differentiation. Differential operators can therefore be applied directly to the basis functions at arbitrary collocation points, and the reduced coefficients can be estimated using the same strong-form, differentiable residual machinery used by PINNs.

Residual-minimizing reduced-order methods determine reduced coefficients by minimizing a governing-equation residual. Examples include least-squares Petrov--Galerkin methods and the physics-informed Karhunen--Lo\`eve expansion (PICKLE) method~\cite{Tartakovsky2020JCP_PICKLE,carlberg2011efficient,carlberg2017galerkin,carlberg2013gnat}. Space--time extensions, including space--time LSPG and dynamic PICKLE (dPICKLE), construct a low-dimensional representation of the complete trajectory and minimize the residual jointly over space and time~\cite{choi2019space,tartakovsky2024physics}. This reduces both spatial and temporal dimensions and avoids evolving a spatial reduced-order model through every full-order time step. In these methods, however, the basis functions, residual operators, and reduced solution are generally defined on the spatial or space--time discretization used to construct the model. lPINN adopts a related residual-minimization strategy but represents the space--time basis functions by coordinate-dependent neural networks. Each basis function depends jointly on the spatial and temporal coordinates, and a single coefficient vector represents the complete trajectory. Because the basis functions are continuous and differentiable through automatic differentiation, the differential-equation residual can be evaluated at arbitrary collocation points, independently of the offline snapshot mesh. This enables cross-resolution evaluation and permits the online collocation set to differ from, and potentially be finer than, the discretization used to generate the training ensemble. For linear differential operators, the online lPINN problem reduces to regularized linear least squares. For nonlinear operators, it remains nonlinear but is restricted to the reduced coefficients and, in inverse problems, the unknown physical parameters. This flexibility requires offline neural-basis training and, as in other snapshot-based methods, depends on the representativeness of the training ensemble.

The continuous trial space connects lPINN to neural fields, also called implicit neural representations. Neural fields map coordinates to physical or geometric quantities and provide continuous, differentiable representations that can be queried independently of the training discretization~\cite{xie2022neural}. Recent continuous reduced-order modeling methods use neural fields to represent solution manifolds and evolve low-dimensional latent states, thereby removing the direct dependence of conventional ROM bases on a particular mesh~\cite{chen2023crom}. Physics-informed conditional neural-field ROMs have further incorporated differential-equation residuals through automatic differentiation and have considered exact treatment of initial and boundary conditions~\cite{kim2024physics}. The lPINN method shares the use of coordinate-based differentiable representations, but its reduced coordinates have a different role. Rather than learning a nonlinear decoder driven by a latent dynamical model, lPINN uses the last hidden-layer features as a linear neural trial basis and determines the instance-specific coefficient vector directly from the governing residual and available observations. This construction preserves a transparent reduced expansion, gives a linear online solve whenever the differential operator and constraints are linear, and naturally supports joint state and parameter estimation without requiring a separately learned latent evolution law or parameter-to-latent map.

A second closely related class comprises extreme-learning and random-feature physics-informed methods. Physics-informed extreme learning machines freeze randomly generated hidden-layer parameters and determine output weights from physics and boundary constraints, leading to rapid linear solves for linear differential equations~\cite{dwivedi2020physics,liu2023bayesian}. Random-feature methods similarly approximate the solution in a prescribed random nonlinear feature space and estimate the coefficients by collocation and least squares~\cite{chen2022random}. These approaches demonstrate the computational advantage of restricting training to linear output weights. Their feature spaces, however, are generally random or problem-independent and may require many features, multiscale constructions, domain decomposition, or careful rescaling to represent a parameterized solution family. In lPINN, the nonlinear features are not random. They are learned offline from an ensemble of solutions and can be trained to reproduce not only the states but also the derivatives or residual quantities required by the governing operator. The learned features therefore form an operator-compatible reduced basis adapted to the solution manifold. Online inference retains the efficiency of output-layer optimization while using a feature space informed by both the solution ensemble and, when desired, the governing physics.

Transfer-learning PINNs also reuse representations across related differential equations. Standard transfer learning initializes a new PINN from parameters trained on a related problem and then fine-tunes some or all network layers. Of particular relevance is the one-shot approach of Desai et al.~\cite{desai2021one}, which freezes the hidden layers of a bundle-trained PINN and adapts the final linear layer for new linear ordinary and partial differential equations. That work establishes that final-layer adaptation can convert a linear forward problem into a linear least-squares solve. lPINN builds on the general efficiency principle of frozen nonlinear features but differs in its objective and scope. Its offline representation is explicitly constructed as a mean-plus-fluctuation reduced basis from numerical solution ensembles, with derivative-matching or residual-based objectives designed to make the basis compatible with the target operator. Its online formulation covers nonlinear forward problems, for which only the reduced coefficients are optimized, and inverse problems, for which coefficients and unknown physical parameters are estimated jointly from physics and sparse data. Consequently, lPINN is better viewed as a physics-corrected neural reduced-basis method than as conventional parameter fine-tuning.

Neural operators take a different route to amortizing PDE solution cost. DeepONet learns an operator through branch and trunk networks, while the Fourier neural operator learns mappings between function spaces using integral-kernel layers parameterized in Fourier space~\cite{lu2021learning,li2021fourier}. Once trained, these methods directly map coefficients, source terms, initial conditions, or boundary data to a solution field and can evaluate new instances without solving a problem-specific optimization. Physics-informed operator-learning methods, including physics-informed DeepONets and physics-informed neural operators, incorporate governing-equation residuals into operator training to reduce data requirements and improve physical consistency~\cite{wang2021learning,li2024physics}. lPINN differs from both purely data-driven and physics-informed neural operators because it does not learn the complete input-to-solution map. The offline data are used to learn a trial space, not a direct predictor of the final state or its reduced coordinates. For each new instance, the coefficients are recomputed by enforcing that instance's governing equation, conditions, and observations. This online physical correction can adapt the prediction to the particular query and incorporate sparse measurements without retraining a global operator. It also avoids constructing an explicit encoding of every possible parameter, source, initial condition, or boundary function. The trade-off is that lPINN requires a small online solve, whereas a neural operator usually provides a single feed-forward prediction.

The lPINN formulation provides several methodological advantages. First, it combines the dimensional reduction and offline--online amortization of reduced-basis methods with a continuous neural representation that can be evaluated and differentiated at arbitrary coordinates. Second, it performs problem-specific residual minimization, unlike surrogates, so the numerical snapshots define the admissible trial space without completely determining the solution at a new query. Third, operator-compatible pretraining can concentrate representational capacity on the states and derivatives that control the online residual, in contrast to random-feature methods. Fourth, the linear dependence on the last-layer coefficients yields a convex online problem for linear operators and substantially reduces the dimension of nonlinear forward and inverse optimization. Fifth, observations can be incorporated directly into the reduced physics objective, allowing state and parameter inference in a unified formulation. These benefits come with the standard limitations of snapshot-based model reduction: the offline ensemble must adequately cover the anticipated solution manifold, and the offline simulation and basis-training costs must be amortized over sufficiently many queries.

The principal contributions of this work are as follows:
\begin{enumerate}
    \item We develop operator-compatible continuous neural basis functions learned from numerical solution ensembles. The coordinate-dependent basis functions are differentiable through automatic differentiation and can be pretrained using state and derivative information or state data augmented by PDE residuals.
    \item We formulate forward prediction as residual-based coefficient inference within the pretrained basis function space. The online problem reduces to regularized linear least squares for linear operators and to reduced-dimensional nonlinear least squares for nonlinear operators, avoiding repeated optimization of the full nonlinear network.
    \item We extend the reduced formulation to inverse problems by jointly estimating the basis function coefficients and unknown physical parameters from the PDE residuals and sparse observations.
    \item We demonstrate compatibility with representation-enhancing architectures, including Fourier features for oscillatory solutions, and investigate derivative-matching and residual-based strategies for constructing operator-compatible neural basis functions.
\end{enumerate}

The numerical experiments consider forward and inverse problems for the advection--diffusion equation, Burgers' equation, and the nonlinear pendulum equation. 
The numerical study evaluates the central hypothesis of this work: that neural representation learning can be separated from problem-specific physics inference so repeated forward and inverse problems can be solved without optimizing the full nonlinear network for each new instance. We therefore use vanilla PINNs as the primary baseline in the following experiments. 
We also provide a comparison with a residual-minimizing ROM \cite{tartakovsky2024physics}.
A systematic comparison with classical projection-based and nonlinear ROMs is outside the scope of the present work and is reserved for future study.
The experiments examine predictive accuracy, parameter estimation, online computational cost, sensitivity to residual and measurement sampling, the effects of basis dimension and regularization, and transfer of the learned continuous representation across evaluation resolutions.

\section{Methods}
\label{sec:surrogate_model}

\subsection{Vanilla Physics-Informed Neural Network (PINN)}

Consider a PDE of the form
\begin{equation}
\label{eq:pde_general}
\mathcal{L}(h(x,t);y)=0, 
\qquad x\in\Omega,\quad t\in (0,T],
\end{equation}
subject to the initial conditions
\begin{equation}
\label{eq:ic_general}
h(x,0)=h_0(x), \quad \partial h/ \partial t(x,0)=  h_1(x)
\qquad x\in\Omega,
\end{equation}
and the boundary condition
\begin{equation}
\label{eq:bc_general}
\mathcal{B}(h(x,t))=g(x,t),
\qquad x\in\Gamma,\quad t\in (0,T],
\end{equation}
where $\mathcal{L}$ is a known differential operator, $\Omega$ is the spatial domain with boundary $\Gamma=\partial\Omega$, $T$ is the time horizon, $h(x,t)$ is the solution, $y$ denotes the physical parameter(s) of the system, $h_0(x)$ and $h_1(x)$ are the zero and first order initial conditions, $\mathcal{B}$ is the boundary operator, and $g(x,t)$ is the prescribed boundary data.

In the PINN framework, the solution $h(x,t)$ is approximated by a neural network $\hat{h}(x,t;\theta)$ with trainable parameters $\theta$. In  forward problems, the physical parameter $y$ is known, and the network parameters are determined by minimizing a residual least-squares loss function:
\begin{align}
\label{eq:pinn_forward_discrete}
\bm{\theta^*}
=
\arg\min_{\bm{\theta}}
\Bigg(
&\frac{1}{N_{\mathrm{res}}}
\sum_{i=1}^{N_{\mathrm{res}}}
\left[
\mathcal{L}\bigl(\hat{h}(x_i^r,t_i^r;\bm{\theta});y\bigr)
\right]^2
\nonumber\\[4pt]
&+
\frac{\lambda_{\mathrm{BC}}}{N_{\mathrm{BC}}}
\sum_{j=1}^{N_{\mathrm{BC}}}
\left[
\mathcal{B}\bigl(\hat{h}(x_j^{\mathrm{BC}},t_j^{\mathrm{BC}};\bm{\theta})\bigr)
-
g(x_j^{\mathrm{BC}},t_j^{\mathrm{BC}})
\right]^2
\nonumber\\[4pt]
&+
\frac{\lambda_{\mathrm{IC},0}}{N_{\mathrm{IC}}}
\sum_{k=1}^{N_{\mathrm{IC}}}
\left[
\hat{h}(x_k^{\mathrm{IC}},0;\bm{\theta})-h_0(x_k^{\mathrm{IC}})
\right]^2
+
\frac{\lambda_{\mathrm{IC},1}}{N_{\mathrm{IC}}}
\sum_{k=1}^{N_{\mathrm{IC}}}
\left[
\frac{\partial \hat{h}(x_k^{\mathrm{IC}},0;\bm{\theta})}{\partial t}-h_1(x_k^{\mathrm{IC}})
\right]^2
+
\lambda\|\boldsymbol{\theta}\|_2^2
\Bigg),
\end{align}
where $\lambda_{\mathrm{BC}}$,  $\lambda_{\mathrm{IC},0}$, and $\lambda_{\mathrm{IC},1}$ are weights associated with the boundary and initial condition terms, respectively; $\lambda$ is the L2 regularization weight; and $N_{\text{res}}$, $N_{\mathrm{BC}}$, and $N_{\mathrm{IC}}$ denote the numbers of collocation points in the interior, boundary, and initial domains, respectively.

In inverse problems, $y$ is treated as an additional trainable variable and is estimated together with $\bm\theta$ using both physics constraints and measurements of $h$,
$\{(X_i^{s},T_i^{s},h_i^{s})\}_{i=1}^{N_{\mathrm{m}}}$, where $X_i^{s}$ and $T_i^{s}$ denote the space and time coordinates of the measurements and $h_i^{s}$ are the measured values. The inverse PINN formulation is given by
\begin{align}\nonumber
(\bm{\theta}^*,y^*) = \arg\min_{\bm{\theta},y}
\Bigg(
&\frac{1}{N_{\text{res}}}
\sum_{i=1}^{N_{\text{res}}}
\left[
\mathcal{L}\bigl(\hat{h}(x_i^r,t_i^r;\bm{\theta});y\bigr)
\right]^2
\nonumber\\[4pt]
&+
\frac{\lambda_{\mathrm{BC}}}{N_{\mathrm{BC}}}
\sum_{j=1}^{N_{\mathrm{BC}}}
\left[
\mathcal{B}\bigl(\hat{h}(x_j^{\mathrm{BC}},t_j^{\mathrm{BC}};\bm{\theta})\bigr)
-
g(x_j^{\mathrm{BC}},t_j^{\mathrm{BC}})
\right]^2
\nonumber\\[4pt] \nonumber
&+
\frac{\lambda_{\mathrm{IC},0}}{N_{\mathrm{IC}}}
\sum_{k=1}^{N_{\mathrm{IC}}}
\left[
\hat{h}(x_k^{\mathrm{IC}},0;\bm{\theta})-h_0(x_k^{\mathrm{IC}})
\right]^2
+\frac{\lambda_{\mathrm{IC},1}}{N_{\mathrm{IC}}}
\sum_{k=1}^{N_{\mathrm{IC}}}
\left[
\frac{\partial \hat{h}(x_k^{\mathrm{IC}};\bm{\theta})}{\partial t}-h_1(x_k^{\mathrm{IC}})
\right]^2
\\ 
&+
\frac{\lambda_{\mathrm{data}}}{N_{\mathrm{m}}}
\sum_{i=1}^{N_{\mathrm{m}}}
\left[
\hat{h}(X_i^{s},T_i^{s};\bm{\theta})-h_i^{s}
\right]^2
+
\lambda\|\boldsymbol{\theta}\|_2^2
\Bigg),
\label{eq:pinn_inverse_continuous}
\end{align}
where the second-to-last term represents the data mismatch between the network prediction and the measurement data. 

The optimization problems in Egs. \eqref{eq:pinn_forward_discrete} and \eqref{eq:pinn_inverse_continuous} are nonlinear least-squares problems in the neural network parameters. As the number of trainable parameters increases, the optimization may become severely ill-conditioned, leading to slow convergence and reduced robustness, especially for stiff and nonlinear PDEs. 

\subsection{Linearized Physics-Informed Neural Network (lPINN)}
\label{sec:lPINN_method}

Our objective is to simplify the training of PINNs by pretraining some of the DNN's parameters using an ensemble of PDE solutions that we refer to as the training dataset:
\[
D_{\mathrm{train}}
=
\left\{
\bigl(g^{(i)},y^{(i)},h_0^{(i)},h_1^{(i)}\bigr)\mapsto \bm{h}^{(i)}
\right\}_{i=1}^{N_{\mathrm{train}}},
\]
where each sample consists of different instances of the PDE parameters, initial, and boundary conditions, and the corresponding PDE solution vector $\bm{h}$ obtained on a mesh with $N_x$ spatial nodes and $N_t$ time steps. 

We start by decomposing the PDE solutions into the sample mean $\bar{h}(x,t)$ and  fluctuations $h'(x,t)$,
\begin{equation}
\label{eq:h_decomp}
h(x,t)=\bar{h}(x,t)+h'(x,t).
\end{equation}
Both $\bar{h}(x,t)$ and $h'(x,t)$ are approximated by DNNs. The fluctuation field is approximated by the DNN
\begin{equation}
\label{eq:fluctuation_nn}
h'(x,t)\approx \hat{h}'(x,t;\bm{\theta})
=
\mathcal{NN}(\bm{z};\bm{\theta}),
\end{equation}
where $\bm{z}=[x,t]^{\mathrm T}$ denotes the input vector and $\bm{\theta}$ is the collection of trainable network parameters. For generality, we consider a fully connected network with $N$ hidden layers,
\begin{equation}
\mathcal{NN}(\bm{z};\bm{\theta})
=
\rho_{N+1}\Bigl(
\rho_N(
\rho_{N-1}(
\cdots
\rho_1(\bm{z})
))
\Bigr),
\end{equation}
where each layer is defined by
\begin{equation}
\bm{z}_{i+1}
=
\rho_i(\bm{W}_i\bm{z}_i+\bm{b}_i),
\qquad i=1,\dots,N,
\end{equation}
with $\bm{z}_1=\bm{z}$. Here, $\bm{W}_i$ and $\bm{b}_i$ are the weights and biases of layer $i$, respectively, and $\rho_i$ is the activation function. In this work, we use the hyperbolic tangent activation function for the hidden layers and the identity function for the output layer. Thus, the full parameter set is
\[
\bm{\theta}=(\bm{W}_{1:N+1},\bm{b}_{1:N+1}).
\]

Because of the choice of activation functions, the network is nonlinear with respect to the hidden-layer parameters
\[
\tilde{\bm{\theta}}=(\bm{W}_{1:N},\bm{b}_{1:N}),
\]
and linear with respect to the output layer parameters $\bm{W}_{N+1}$ and $\bm{b}_{N+1}$. We denote the size of the last hidden layer by $N_\eta$ and note that the output layer has a single output for a scalar state variable. Then, $\bm{W}_{N+1}$ is the $N_\eta$-dimensional vector and $\bm{b}_{N+1}$ is a scalar. 
Because we use the DNN to model fluctuations around the mean, we set $\bm{b}_{N+1}=0$. 
Finally, we can express fluctuations as
\begin{equation}
\label{eq:linear_last_layer}
\hat{h}'(\bm{z};\bm{\theta})
=
\hat{h}'(\bm{z};\tilde{\bm{\theta}},\bm{W}_{N+1})
=
\bm{W}_{N+1} \cdot \bm\psi(\bm{z};\tilde{\bm{\theta}}),
\end{equation}
where $\bm\psi = [\psi_1,...,\psi_{N_\eta}]^T$ and $\psi_i(\bm{z};\tilde{\bm{\theta}})$ is the output of $i$th neuron of the last hidden layer. 

Eq \eqref{eq:linear_last_layer} provides an expansion of $h'(x,t)$ in terms of the product of space-time-dependent parametric basis functions $\psi_i(x,t;\tilde{\bm{\theta}})$ and coefficients $\bm{W}_{N+1}$. Because $\psi_i(x,t;\tilde{\bm{\theta}})$ are defined by the DNN layers, we refer to them as neural basis functions.

As in finite element, POD-ROM, and other reduced-order methods, we assume that the solution manifold can be approximated by a set of basis functions that are independent of the PDE parameters, initial conditions, and boundary conditions. The effect of PDE parameters, initial, and boundary conditions on the solution is captured by $\bm{W}_{N+1}$.

In standard PINN, $\tilde{\bm{\theta}}$ and $\bm{W}_{N+1}$ are learned jointly by minimizing the PDE residuals, which results in a highly nonlinear least-squares minimization problem.  

The key idea of lPINN is to learn the basis functions, i.e., to estimate 
$\tilde{\bm\theta}$, from a training dataset 
$\{ y^{(i)} \rightarrow \bm{h}^{(i)} \}_{i=1}^{N_{\mathrm{train}}}$ 
during an ``offline'' stage, and only estimate $\bm{W}_{N+1}$ ``online'' by minimizing the PDE residuals. In the training dataset, $y^{(i)}$ are the samples of the parameters from the desired range and $\bm{h}^{(i)}$ are the corresponding solutions obtained (in general) numerically on the mesh with $N_e$ elements or grid points.  

The mean field is approximated by another DNN that is also trained offline:
\begin{equation}
\label{eq:mean_nn}
\bar{h}(x,t)\approx \hat{\bar{h}}(x,t;\bm{\gamma})
=
\mathcal{NN}(\bm{z};\bm{\gamma}),
\end{equation}
where $\bm{\gamma}$ is a collection of trainable parameters.
The DNN $\hat{\bar{h}}(x,t;\bm{\gamma})$ is trained using the sample mean values:
\begin{equation}
\label{eq:h_mean}
\bar{\bm{h}}
\approx
\frac{1}{N_{\mathrm{train}}}
\sum_{i=1}^{N_{\mathrm{train}}}
\bm{h}^{(i)}.
\end{equation}
We would like $\hat{\bar{h}}(\bm{x},t;\bm{\gamma})$ and $\hat{h}'(\bm{x},t;\bm{\theta})$ to provide an accurate approximation of not just ${h}(x,t;\bm{\gamma})$ but also all its derivatives present in the governing equation. Assume, for generality, that $\mathcal{L}$ contains first and second space and time derivatives. Then, we supplement the $\bar{\bm{h}}$ data with its derivatives \(\partial_t \bar{\bm{h}}\), \(\partial_x \bar{\bm{h}}\), \(\partial_{xx} \bar{\bm{h}}\), and \(\partial_{tt} \bar{\bm{h}}\). These derivatives can be computed numerically from $\bar{\bm{h}}$ on the same space-time mesh that was used to compute $\bar{\bm{h}}$.

Then, $\bm{\gamma}$ can be obtained by solving an optimization problem, which minimizes the mismatch between $\hat{\bar{h}}(\bm{z};\bm{\gamma})$, $\bar{\bm{h}}$, and their derivatives:  
\begin{equation}
\label{eq:gamma_opt}
\begin{aligned}
\bm{\gamma}^{*}
=
\arg\min_{\bm{\gamma}}
\Big[
&\|\hat{\bar{\bm{h}}}(\bm{\gamma})-\bar{\bm{h}}\|_2^2
+\lambda_1
\|\partial_t \hat{\bar{\bm{h}}}(\bm{\gamma})-\partial_t \bar{\bm{h}}\|_2^2
\\
&+\lambda_2
\|\partial_x \hat{\bar{\bm{h}}}(\bm{\gamma})-\partial_x \bar{\bm{h}}\|_2^2
+\lambda_3
\|\partial_{xx} \hat{\bar{\bm{h}}}(\bm{\gamma})-\partial_{xx} \bar{\bm{h}}\|_2^2
\\
&+\lambda_4
\|\partial_{tt} \hat{\bar{\bm{h}}}(\bm{\gamma})-\partial_{tt} \bar{\bm{h}}\|_2^2
+
\lambda\|\boldsymbol{\gamma}\|_2^2
\Big].
\end{aligned}
\end{equation}
where \(\hat{\bar{\bm{h}}}(\bm{\gamma})\),  
\(\partial_t \hat{\bar{\bm{h}}}(\bm{\gamma})\), 
\(\partial_x \hat{\bar{\bm{h}}}(\bm{\gamma})\),
\(\partial_{xx} \hat{\bar{\bm{h}}}(\bm{\gamma})\), and \(\partial_{tt} \hat{\bar{\bm{h}}}(\bm{\gamma})\)
are the vectors of the $\hat{\bar{h}}(x,t;\bm{\gamma})$ values evaluated on the space-time mesh, as well as the derivatives of $\hat{\bar{h}}(x,t;\bm{\gamma})$ evaluated on the same mesh. The weights  $\lambda_1$, $\lambda_2$, $\lambda_3$, $\lambda_4$ and $\lambda$ are used to balance different terms and are selected through a grid search. 

Similarly, $\tilde{\bm{\theta}}$ is estimated from the dataset $\{ \bm{h}^{\prime(i)},
\partial_t \bm{h}^{\prime(i)},
\partial_{tt} \bm{h}^{\prime(i)},
\partial_x \bm{h}^{\prime(i)},
\partial_{xx} \bm{h}^{\prime(i)}
\}_{i=1}^{N_{\mathrm{train}}}$,
where $\bm{h}^{\prime(i)} = \bm{h}^{(i)} - \overline{\bm{h}}$.

The coefficients $\{\bm{W}_{N+1}^{(i)}\}_{i=1}^{N_{\mathrm{train}}}$ are found by solving the minimization problem:
\begin{equation}
\label{eq:theta_pretrain}
\begin{aligned}
\tilde{\bm{\theta}}^{*}, \{\bm{W}_{N+1}^{(i)}\}_{i=1}^{N_{\mathrm{train}}}
=
\arg\min_{\substack{
\tilde{\bm{\theta}},\\
\{\bm{W}_{N+1}^{(i)}\}_{i=1}^{N_{\mathrm{train}}}
}}
\sum_{i=1}^{N_{\mathrm{train}}}
\Big[
&
\|\hat{\bm{h}}^{\prime(i)}(\tilde{\bm{\theta}},
\bm{W}_{N+1}^{(i)} )-\bm{h}^{\prime(i)}\|_2^2
+\lambda_1
\|\partial_t \hat{\bm{h}}^{\prime(i)}(\tilde{\bm{\theta}},
\bm{W}_{N+1}^{(i)} ) -\partial_t \bm{h}^{\prime(i)}\|_2^2
\\
&+\lambda_2
\|\partial_x\hat{\bm{h}}^{\prime(i)}(\tilde{\bm{\theta}},
\bm{W}_{N+1}^{(i)} ) - \partial_x \bm{h}^{\prime(i)}\|_2^2
+\lambda_3
\|\partial_{xx}\hat{\bm{h}}^{\prime(i)}(\tilde{\bm{\theta}},
\bm{W}_{N+1}^{(i)} )-\partial_{xx} \bm{h}^{\prime(i)}\|_2^2
\\
&+\lambda_4
\|\partial_{tt}\hat{\bm{h}}^{\prime(i)}(\tilde{\bm{\theta}},
\bm{W}_{N+1}^{(i)} )-\partial_{tt} \bm{h}^{\prime(i)}\|_2^2
\Big]
+
\lambda
\left(
\|\tilde{\boldsymbol{\theta}}\|_2^2
+
\sum_{i=1}^{N_{\mathrm{train}}}
\|\mathbf{W}_{N+1}^{(i)}\|_2^2
\right)
\end{aligned}
\end{equation}
where $\hat{\bm{h}}^{\prime(i)}(\tilde{\bm{\theta}},
\bm{W}_{N+1}^{(i)} )$ is the vector of $\hat{h}'\bigl(\bm{z};
\tilde{\bm{\theta}},
\bm{W}^{(i)}_{N+1}
\bigr)$ predictions evaluated on the space-time mesh $N_x \times N_t$ and $\partial_t\hat{\bm{h}}^{\prime(i)}$, $\partial_{tt} \hat{\bm{h}}^{\prime(i)}$, $\partial_x\hat{\bm{h}}^{\prime(i)}$, and $\partial_{xx} \hat{\bm{h}}^{\prime(i)}$ are the vectors of the derivatives of $\hat{h}'\bigl(\bm{z};
\tilde{\bm{\theta}},
\bm{W}^{(i)}_{N+1}
\bigr)$ evaluated on the same mesh. The estimation of $\tilde{\bm\theta}$ from Eq \eqref{eq:theta_pretrain} also requires estimating $\{ \bm{W}^{(i)}_{N+1} \}_{i=1}^{N_{\mathrm{train}}}$, the parameters $\bm{W}_{N+1}$ specific to each sample in the training dataset, even though these parameters are not part of the pretrained DNN network.

Evaluating derivatives of the solution incurs additional computational cost. Although this cost is substantially smaller than the cost of generating the solution ensemble, it can still increase offline training time and storage requirements and should therefore be avoided when possible. We find that for some PDE types, $\bm{\gamma}$ and $\tilde{\bm{\theta}}$ can be evaluated jointly using the samples $\{ \bm{h}^{(i)} \}_{i=1}^{N_{\mathrm{train}}}$ and PDE residual constraints as
\begin{equation}
\begin{aligned}
\bm{\gamma}^{*},
\tilde{\bm{\theta}}^{*},
\left\{
\bm{W}_{N+1}^{*(i)}
\right\}_{i=1}^{N_{\mathrm{train}}}
=
&
\arg\min_{\substack{
\bm{\gamma},
\tilde{\bm{\theta}},\\
\{\bm{W}_{N+1}^{(i)}\}_{i=1}^{N_{\mathrm{train}}}
}}
\sum_{i=1}^{N_{\mathrm{train}}}
\Big[
\|\hat{\bar{\bm{h}}}(\bm{\gamma})-\bar{\bm{h}}\|_2^2
+\|\hat{\bm{h}}^{\prime(i)}(\tilde{\bm{\theta}},
\bm{W}_{N+1}^{(i)} )-\bm{h}^{\prime(i)}\|_2^2
\Big]
\\
&
+\lambda_f\|\bm{\mathcal{R}}^{(i)} (\bm\gamma,\tilde{\bm{\theta}},
\bm{W}_{N+1}^{(i)} ) \|_2^2
+\lambda_{\mathrm{IC},0}\|\bm{\mathcal{R}}^{(i)}_{\mathrm{IC},0} (\bm\gamma,\tilde{\bm{\theta}},
\bm{W}_{N+1}^{(i)} )\|_2^2
+
\lambda_{\mathrm{IC},1}\|\bm{\mathcal{R}}^{(i)}_{\mathrm{IC},1} (\bm\gamma,\tilde{\bm{\theta}},
\bm{W}_{N+1}^{(i)} )\|_2^2
\\
&+
\lambda_{\mathrm{BC}}\|\bm{\mathcal{R}}^{(i)}_{\mathrm{BC}}(\bm\gamma,\tilde{\bm{\theta}},
\bm{W}_{N+1}^{(i)} )\|_2^2
\Big]
+
\lambda
\left(
\|\boldsymbol{\gamma}\|_2^2
+
\|\tilde{\boldsymbol{\theta}}\|_2^2
+
\sum_{i=1}^{N_{\mathrm{train}}}
\|\mathbf{W}_{N+1}^{(i)}\|_2^2
\right)
\label{eq:joint_train_physics}
\end{aligned}
\end{equation}
where the residual operators are defined as
\begin{equation}
\label{eq:R_pde}
\mathcal{R}(x,t;y;\bm\gamma,\tilde{\bm\theta},\bm{W}_{N+1})
=
\mathcal{L}\Big(
\hat{\bar{h}}(x,t;\bm{\gamma})
+
\hat{h}'(x,t;\tilde{\bm{\theta}},\bm{W}_{N+1});
\, y
\Big),
\end{equation}
\begin{equation}
\label{eq:R_bc}
\mathcal{R}_{\mathrm{BC}}(x,t;g;\bm\gamma,\tilde{\bm\theta},\bm{W}_{N+1})
=
\mathcal{B}\Big(
\hat{\bar{h}}(x,t;\bm{\gamma})
+
\hat{h}'(x,t;\tilde{\bm{\theta}},\bm{W}_{N+1})
\Big)
-
g(x,t),
\end{equation}
\begin{equation}
\label{eq:R_ic}
\mathcal{R}_{\mathrm{IC},0}(x;h_0;\bm\gamma,\tilde{\bm\theta},\bm{W}_{N+1})
=
\hat{\bar{h}}(x,0;\bm{\gamma})
+
\hat{h}'(x,0;\tilde{\bm{\theta}},\bm{W}_{N+1})
-
h_0(x),
\end{equation}
and
\begin{equation}
\label{eq:R_ic1}
\mathcal{R}_{\mathrm{IC},1}(x;h_0;\bm\gamma,\tilde{\bm\theta},\bm{W}_{N+1})
=
\frac{\partial \hat{\bar{h}}(x,0;\bm{\gamma})}{\partial t}
+
\frac{\partial \hat{h}'(x,0;\tilde{\bm{\theta}},\bm{W}_{N+1})}{\partial t}
-
h_1(x).
\end{equation}
Here, $\bm{\mathcal{R}}^{(i)}$ is the vector with components  ${\mathcal{R}}(x,t;y^{(i)};\bm\gamma,\tilde{\bm\theta},\bm{W}^{(i)}_{N+1})$ evaluated at randomly sampled residual collocation points on the space-time domain; 
${\mathcal{R}}^{(i)}_{\mathrm{IC},0}$ and ${\mathcal{R}}^{(i)}_{\mathrm{IC},1}$ are the vector with components
${\mathcal{R}}_{\mathrm{IC},0}(x;h_0^{(i)};\bm\gamma,\tilde{\bm\theta},\bm{W}_{N+1}^{(i)})$ and ${\mathcal{R}}_{\mathrm{IC},1}(x;h_1^{(i)};\bm\gamma,\tilde{\bm\theta},\bm{W}_{N+1}^{(i)})$, respectively, evaluated on the space domain at $t=0$;
$\bm{\mathcal{R}}^{(i)}_{\mathrm{BC}}$ is the vector with  components ${\mathcal{R}}_{\mathrm{BC}}(x,t;g^{(i)};\bm\gamma,\tilde{\bm\theta},\bm{W}_{N+1}^{(i)})$ evaluated on the time domain on the boundary; and 
$y^{(i)}$, $g^{(i)}$, $h_0^{(i)}$, $h_1^{(i)}$ are components of the $D_{train}$ training dataset.

After the offline pretraining stage, the PDE solution for any $y$, $g$, and $h_0$ is obtained by estimating $\mathbf{W}_{N+1}$ from the minimization problem

\begin{eqnarray}\label{eq:RLS_final}
\mathbf{W}_{N+1}^{*}
=
\arg\min_{\mathbf{W}_{N+1}}
\Big[&&
||\bm{\mathcal{R}}(y;\gamma^{*},\tilde{\bm{\theta}}^{*},\mathbf{W}_{N+1})\|_2^2
+
\lambda_{\mathrm{BC}}\|\bm{\mathcal{R}}_{\mathrm{BC}}(g;\gamma^{*},\tilde{\bm{\theta}}^{*},\mathbf{W}_{N+1})\|_2^2 \\ \nonumber
&+&
\lambda_{\mathrm{IC},0}\|\bm{\mathcal{R}}_{\mathrm{IC},0}(h_0;\gamma^{*},\tilde{\bm{\theta}}^{*},\mathbf{W}_{N+1})\|_2^2
+
\lambda_{\mathrm{IC},1}\|\bm{\mathcal{R}}_{\mathrm{IC},1}(h_1;\gamma^{*},\tilde{\bm{\theta}}^{*},\mathbf{W}_{N+1})\|_2^2
+
\lambda\|\mathbf{W}_{N+1}\|_2^2
\Big],
\end{eqnarray}
where $\gamma^{*}$ and $\tilde{\bm{\theta}}^{*}$ are obtained at the offline pretraining. The basis functions approximately satisfy the initial and boundary conditions. However, our numerical results show that adding the initial and boundary condition penalty terms, i.e.,  setting $\lambda_{\mathrm{BC}}>0$, $\lambda_{\mathrm{IC},0}>0$, $\lambda_{\mathrm{IC},1}>0$ in general, yields more accurate solutions.

Alternatively, we can use the Galerkin projection method to estimate $\mathbf{W}_{N+1}$ as a solution of the equation
\begin{equation}
\bm{\Psi}_r^{T}
\bm{\mathcal{R}}
\left(
y;\gamma^{*},\bm{\tilde{\theta}^{*}},\mathbf{W}_{N+1}^{*}
\right)
=
\mathbf{0},
\label{eq:galerkin_online}
\end{equation}
where $\bm{\Psi}_r\in \mathbb{R}^{N_{\mathrm{res}}\times N_\eta}$ is the matrix of pretrained neural basis functions evaluated at the residual points, whose $i$th row is $\bm{\psi}(z_i;\tilde{\theta}^{*})^{T}$. 
In our numerical experiments, the Galerkin formulation produced errors similar to those obtained by residual least-square lPINN formulation in Eq \eqref{eq:RLS_final}. We therefore use the residual least-squares lPINN formulation in the results reported below.

In lPINN, the term ``linearized'' refers to the linear dependence of the functional representation on $\bm{W}_{N+1}$, the parameter vector estimated in the PINN-like online stage. It does not imply linearizing the governing differential operator. 
If $\mathcal{L}$ and $\mathcal{B}$ are linear operators, then Eq. \eqref{eq:RLS_final} 
is a regularized linear least-squares problem. It admits a unique solution when the augmented least-squares matrix has full column rank or when strictly positive quadratic regularization makes the objective strictly convex.
For nonlinear operators, the reduced problem remains nonlinear, but it is still low-dimensional and can be solved iteratively using standard nonlinear least-squares methods. 

\begin{figure}[H]
    \centering
    \includegraphics[width=0.95\textwidth]{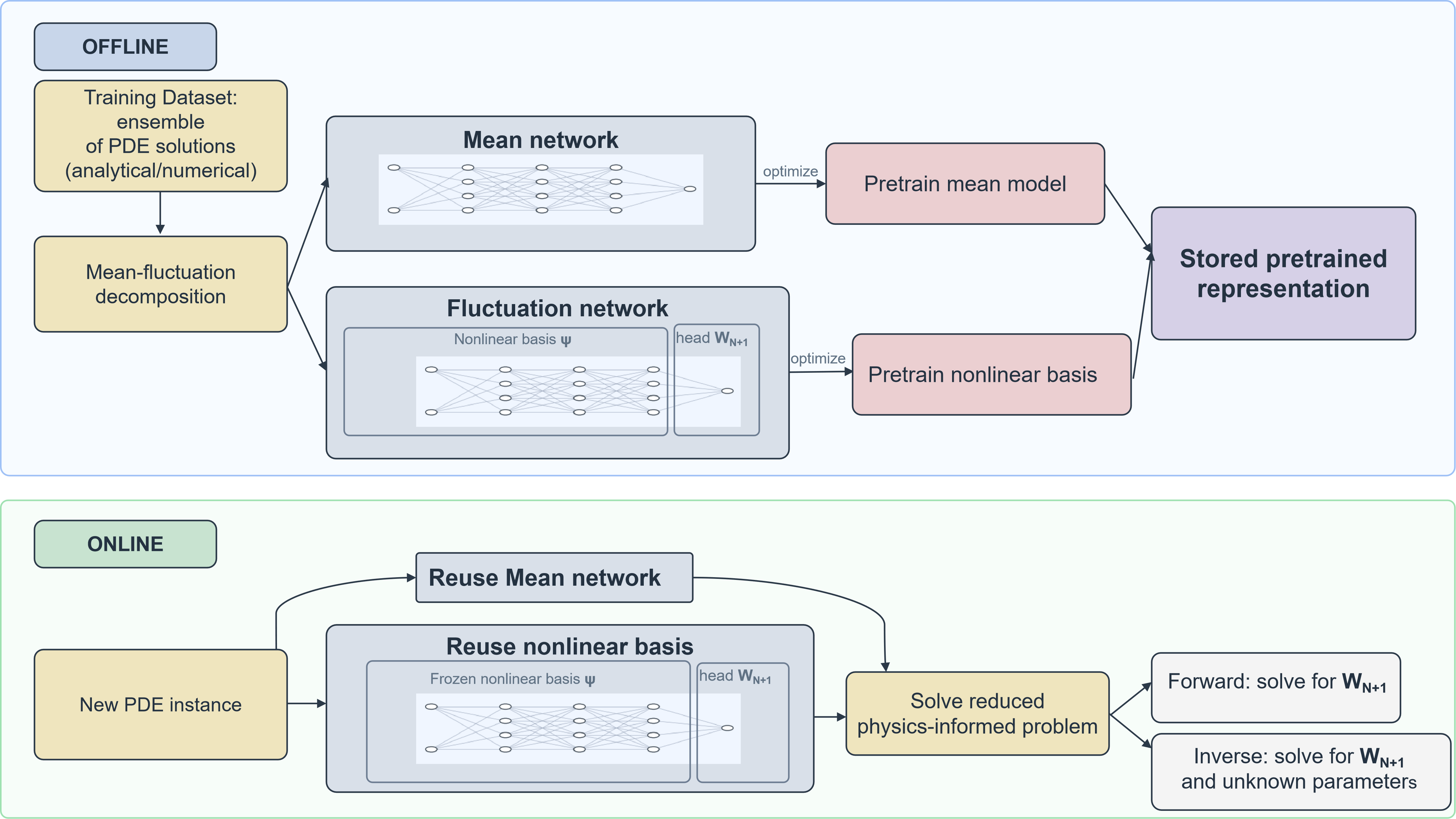}
    \caption{The schematic representation of the lPINN method. All parameters in the mean DNN and the parameters associated with the fluctuation DNN layers are estimated offline using the training dataset. The remaining fluctuation-DNN parameters in the last layer are determined for the specified parameters using the physics-informed loss.
}
    \label{fig:lPINN_flowchart}
\end{figure}

The inverse PDE problem is formulated in lPINN similarly to the PINN method. Consider an inverse problem in which one or more of the parameters $y$, $g$, $h_0$, and $h_1$ are unknown and must be inferred from the observations $\bm{u}^*=[u_j^*]_{j=1}^{N_{\mathrm{m}}}$ of $h$ at locations $\bm{x}_d=[x_j]_{j=1}^{N_{\mathrm{m}}}$ and time instances  $\bm{t}_d=[t_j]_{j=1}^{N_{\mathrm{m}}}$. The parameters can be estimated jointly with $\bm{W}_{N+1}$ at the online stage as the solution of the minimization problem:

\begin{equation}
\begin{aligned}
\left( \bm{W}_{N+1}^{*},\, y^{*}, g^*,h_0^* \right)
=
\arg\min_{\bm{W}_{N+1},\, y, g, h_0}
\Big[
&\|\bm{\mathcal{R}}(y;\bm{W}_{N+1})\|_2^2
+\lambda_{\mathrm{BC}}\|\bm{\mathcal{R}}_{BC}(g;\bm{W}_{N+1})\|_2^2
+\lambda_{\mathrm{IC},0}\|\bm{\mathcal{R}}_{\mathrm{IC},0}(h_0,\bm{W}_{N+1},y)\|_2^2
\\
&
+\lambda_{\mathrm{IC},1}\|\bm{\mathcal{R}}_{\mathrm{IC},1}(h_1,\bm{W}_{N+1},y)\|_2^2
+
\lambda\|\mathbf{W}_{N+1}\|_2^2
+\lambda_{\mathrm{data}}
\left\|
\hat{\bm{h}}_d(\bm{W}_{N+1})-\bm{h}^{*}
\right\|_2^2
\Big],
\end{aligned}
\label{eq:lPINN_inverse}
\end{equation}
where 
L2 regularization is used on the $\mathbf{W}_{N+1}$ parameters, and 
$\bm{h}_d (\bm{W}_{N+1})$ in the data term is the vector with components $\hat{\bar{h}}(x_i,t_i;\bm{\gamma}^*)
+
\hat{h}'(x_i,t_i;\tilde{\bm{\theta}}^*,\bm{W}_{N+1})$. 
Because the basis functions approximately satisfy the initial and boundary conditions, we find that in this loss function, $\lambda_{\mathrm{BC}}$, $\lambda_{\mathrm{IC},0}$ and $\lambda_{\mathrm{IC},1}$ can be set to zero without incurring significant errors in the inverse solution.
The weighting coefficients in the loss function, including $\lambda$ and $\lambda_{\mathrm{data}}$, can be estimated as in PINN. In this work, we use a grid search method to select the optimal combination of the weighting coefficients. The schematic representation of the lPINN method for forward and inverse problems is shown in Figure \ref{fig:lPINN_flowchart}.

A key consequence of the neural representation is that the mesh used to generate the offline solution ensemble and the collocation sets used during offline and online training are decoupled.
Once the mean function $\hat{\overline{\bm{h}}}(\bm{z},\bm\gamma^*)$ and basis functions $\bm\psi(\bm{z};\tilde{\bm{\theta}}^*)$ have been learned, both 
$
\hat{h}(\bm{z};\bm{W}_{N+1}) = \hat{\overline{\bm{h}}}(\bm{z},\bm\gamma^*) + 
 \bm\psi(\bm{z};\tilde{\bm{\theta}}^*) \cdot \bm{W}_{N+1}
$
together with the derivatives required by the differential operator, can be evaluated at arbitrary coordinates $\bm{z}$ through automatic differentiation. Consequently, a basis learned from a solution ensemble generated on a coarse mesh can be used to represent and evaluate the solution at a finer resolution without retraining. The results presented below indicate that this cross-resolution accuracy does not require a large number of residual points during online inference because only a relatively small number of reduced coefficients is estimated. Instead, the attainable accuracy is governed primarily by the quality of the pretrained mean and basis functions. Within the range investigated here, this representation accuracy can be improved by increasing the number of solution samples in the offline training ensemble. Increasing the number of residual points during residual-based offline pretraining can also improve the learned representation, although we observed a smaller improvement than that obtained by increasing the number of training samples. We define lPINN superresolution as the ability of an lPINN trained on a coarse-resolution ensemble to produce predictions that are more accurate than numerical solutions computed on finer meshes when both are evaluated against the same high-fidelity reference solution.

\section{Numerical Examples}\label{sec:numerical_results}

Unless otherwise stated, all reported solution and parameter errors are computed from the arithmetic-mean prediction or parameter estimate over 10 independent runs, and the reported computational times denote the average online training time per run. 

In the comparison studies, the PINN network and the mean and fluctuation networks used in lPINN have the same size and architecture, including the same number of hidden layers, the same width in each hidden layer, and the same activation functions. Thus, each individual network used in lPINN is matched in size and architecture to the PINN network. The exception is a comparison in Section \ref{sec:pendulum} where PINN and lPINN use different Fourier frequency scales.

\subsection{Advection Diffusion Equation}

We first evaluate the proposed lPINN framework on the one-dimensional advection--diffusion equation
\begin{equation}
\label{eq:ADE}
\frac{\partial u}{\partial t}
+
V\frac{\partial u}{\partial x}
=
D\frac{\partial^2 u}{\partial x^2},
\end{equation}
subject to the initial condition
\begin{equation}
\label{eq:ADE_ic}
u(x,0)=0,
\end{equation}
and the boundary conditions
\begin{equation}
\label{eq:ADE_bc_left}
u(0,t)=u_0,
\end{equation}
where $u_0=1$,
and
\begin{equation}
\label{eq:ADE_bc_right}
\left.\frac{\partial u}{\partial x}\right|_{x=L}=0.
\end{equation}

In the forward problem, the goal is to predict the solution for $V\in[V_{\min},V_{\max}]$ and $D\in[D_{\min},D_{\max}]$. In the inverse problem, the goal is to infer $V$ and $D$ from sparse measurements of $u(x,t)$. 

For $t\ll\frac{L}{V}$, the solution of this equation can be approximated by the Ogata--Banks formula,
\begin{equation}
\label{eq:OgataBanks}
u(x,t)
=
\frac{u_0}{2}
\left[
\operatorname{erfc}\left(\frac{x-Vt}{2\sqrt{Dt}}\right)
+
\exp\left(\frac{Vx}{D}\right)
\operatorname{erfc}\left(\frac{x+Vt}{2\sqrt{Dt}}\right)
\right].
\end{equation}

To construct the training dataset, we treat $V$ and $D$ as independent uniformly distributed random variables defined on the intervals $V\in[V_{\min},V_{\max}]$ and $D\in[D_{\min},D_{\max}]$, respectively. Next, we generate $N_{\mathrm{train}}$ samples of $V$ and $D$ from these distributions. For each sample of $V$ and $D$,
the Ogata-Banks solution is evaluated on a uniform $N_x\times N_t$ mesh with nodes $(x_i,t_k)$ excluding $t=0$, where $x_1=0$, $x_{N_x}=L$, $t_1=T/(N_t-1)$, and $t_{N_t}=T$. 
The sample mean is then estimated as
\begin{equation}
\bar{u}(x_i,t_k)
=
\frac{1}{N_{\mathrm{train}}}
\sum_{n=1}^{N_{\mathrm{train}}}
u^{(n)}(x_i,t_k),
\end{equation}
and the fluctuation field is defined by
\begin{equation}
u'(x,t)=u(x,t)-\bar{u}(x,t).
\end{equation}

In all ADE experiments, we set $L=86$, $T=200$. The training dataset contains $N_{\mathrm{train}}=500$ realizations generated from
\[
V\in[0.6V^{*},\,1.4V^{*}],
\qquad
D\in[0.6D^{*},\,1.4D^{*}],
\]
where $V^{*}=0.213$ and $D^{*} = 0.061$ yielding the Peclet number 
$\mathrm{Pe}=\frac{V^* L}{D^*}=300$. 
The trained model is tested at
\[
V^{ref}=1.2V^{*},
\qquad
D^{ref}=0.7D^{*},
\]
which corresponds to $\mathrm{Pe}=514$. The test case is intentionally selected away from the center of the training ranges. Also, the ADE solutions for $\mathrm{Pe}\gg 1$ are challenging to obtain with standard numerical methods. Therefore, the problem considered provides a rigorous test for lPINN.

The accuracy of lPINN and other methods for ADE is estimated using the dimensionless relative root-mean-square errors (rRMSEs) for the predicted solution $\mathbf{u}$ and, for the inverse problem, the estimated parameters $V$ and $D$. The relative errors are defined as
\[
\mathrm{rRMSE}_u
=
\frac{
\left\|
\mathbf{u}-\mathbf{u}^{\mathrm{ref}}
\right\|_2
}{
u_0\sqrt{N_{\mathrm{u}}}
},
\qquad
\mathrm{rRMSE}_V
=
\frac{
\left|V-V^{\mathrm{ref}}\right|
}{
\left|V^{\mathrm{ref}}\right|
},
\qquad
\mathrm{rRMSE}_D
=
\frac{
\left|D-D^{\mathrm{ref}}\right|
}{
\left|D^{\mathrm{ref}}\right|
},
\]
where $N_u$ is the number of grid points where the numerical solutions are compared to the reference solution $\bm{u}^{\text{ref}}$ (generally, $N_u=N_g$), where $N_g$ denotes the number of grid points on the mesh used to evaluate the solution, and $u_0=1$ is the maximum value of the reference solution achieved at $x=0$. The dimensionless rRMSEs for the PDE residual $\bm{\mathcal{R}}$, initial condition, and boundary conditions are defined as
\[
\mathrm{rRMSE}_R
=
\frac{L}{u_0 V}
\frac{\left\|\bm{\mathcal{R}}\right\|_2}{\sqrt{N_{\mathrm{res}}}},
\qquad
\mathrm{rRMSE}_{\mathrm{IC}}
=
\frac{
\left\|
\mathbf{u}_{\mathrm{IC}}-\mathbf{u}_{\mathrm{IC}}^{\mathrm{ref}}
\right\|_2
}{
u_0\sqrt{N_{\mathrm{IC}}}
},
\]
and
\[
\mathrm{rRMSE}_{BC}
=
\frac{
\left\|
\mathbf{u}(x=0)-\mathbf{u}^{\mathrm{ref}}(x=0)
\right\|_2
}{
u_0\sqrt{N_{\mathrm{BC},0}}
}
+
\frac{L}{
u_0
}
\frac{
\left\|
\partial_x\bm{u}(x=L)
\right\|_2
}{
\sqrt{N_{\mathrm{BC},L}}
},
\]
where $N_{\mathrm{IC}}$ is the number of residual points at $t=0$ (excluding the point $(x,t)=(0,0)$), $N_{\mathrm{BC},0}$ is the number of residual points over the time domain on the $x=0$ boundary, and $N_{\mathrm{BC},L}$ is the number of residual points on the $x=L$ boundary.

In lPINN, the mean and fluctuation networks are pretrained using the Adam optimizer. The online training for the forward ADE problem, i.e., solving the minimization problem \eqref{eq:RLS_final}, is a linear least-squares problem which we solve using the linear least-squares solver numpy.linalg.lstsq. 
The inverse ADE minimization problem is nonlinear and is solved using the L-BFGS method. The number of residual points is set to $N_{\mathrm{res}}=10 N_\eta$. Inverse solutions are obtained with $N_{\mathrm{m}}=40$ measurements of $u$.

We compared the derivative-matching (Eqs.~\eqref{eq:gamma_opt} and \eqref{eq:theta_pretrain}) and
residual-based (Eq.~\eqref{eq:joint_train_physics}) offline pretraining methods. 
For residual-based offline pretraining, we split the $N_x \times (N_t-1)=870$ grid points into 10 batches at each epoch. For each batch, solution values of $N_{\mathrm{train}}=500$ realizations at 87 grid points are used for training, together with $N_{\mathrm{res}}=435$ residual points, $N_{\mathrm{IC}}=435$ initial condition collocation points, and $N_{\mathrm{BC}}=435$ boundary condition collocation points randomly resampled for each batch.
For derivative-matching offline pretraining, the same 870 grid points are split into 10 batches, and for each batch, the solution and derivative values at the corresponding 87 grid points are used for training, together with $N_{\mathrm{IC}}=435$ and $N_{\mathrm{BC}} = N_{\mathrm{BC},0}+N_{\mathrm{BC},L}=435$ randomly resampled initial condition and boundary condition collocation points.
The number of neural basis functions is set to $N_\eta = 50$. The performance of the models in terms of rRMSE of forward and inverse solutions is given in Table \ref{tab:ADE_summary}.  
The lPINN forward solution error is a combination of the approximation error and residual least-squares error. The former measures how well the pretrained DNNs approximate the PDE solution. The latter measures how accurately the residual least-squares formulation estimates $\bm{W}_{N+1}$. To compute the approximation error, we estimate $\bm{W}_{N+1}$ by least-squares fitting the pretrained DNNs to the reference PDE solution. The residual-based pretraining gives a slightly smaller approximation error, $3.94\times 10^{-4}$, compared with $4.44\times 10^{-4}$ for the derivative-matching pretraining. Another related metric is the value of the PDE residual yielded by the pretrained DNN with $\bm{W}_{N+1}$ obtained from the least-squares fitting to the reference PDE solution. The RMSE of the residuals (relative to the zero-residual value produced by the reference solution) is smaller for residual-based pretraining ($3.60\times 10^{-2}$) than in the derivative-matching pretraining ($4.36\times 10^{-2})$. 
However, the residual-based pretraining gives a much smaller lPINN forward solution error, $1.45\times 10^{-3}$ versus $9.38\times 10^{-3}$ in the derivative-matching pretraining. In the inverse problem, residual-based pretraining gives better accuracy for the solution $u$ and velocity $V$, while derivative-matching pretraining gives a smaller error for the diffusion coefficient $D$. Since we seek the best overall performance and, in particular, accurate forward prediction, we conclude that residual-based pretraining is the better strategy for the considered ADE problem, and use it in the remaining numerical examples of this section.

\begin{table}[H]
\centering
\caption{ADE comparison between residual-based pretraining and derivative-matching pretraining. The smaller value in each row is shown in bold. }
\begin{tabular}{lcc}
\hline
Case / Metric & Residual-based & Derivative-matching \\
\hline
Approximation error: rRMSE& $\mathbf{3.94\times10^{-4}}$& $4.44\times10^{-4}$\\
Residual: rRMSE& $\mathbf{3.60\times10^{-2}}$& $4.36\times10^{-2}$\\
Forward: rRMSE$_u$ & $\mathbf{1.45\times10^{-3}}$& $9.38\times10^{-3}$\\
Inverse: rRMSE$_u$ & $\mathbf{3.90\times10^{-3}}$& $2.01\times10^{-2}$\\
Inverse: rRMSE$_V$ & $\mathbf{4.13\times10^{-3}}$& $1.22\times10^{-2}$\\
Inverse: rRMSE$_D$ & $1.07\times10^{-1}$& $\mathbf{6.23\times10^{-2}}$\\
\hline
\end{tabular}
\label{tab:ADE_summary}
\end{table}

\begin{figure}[H]
    \centering
    \includegraphics[width=1\linewidth]{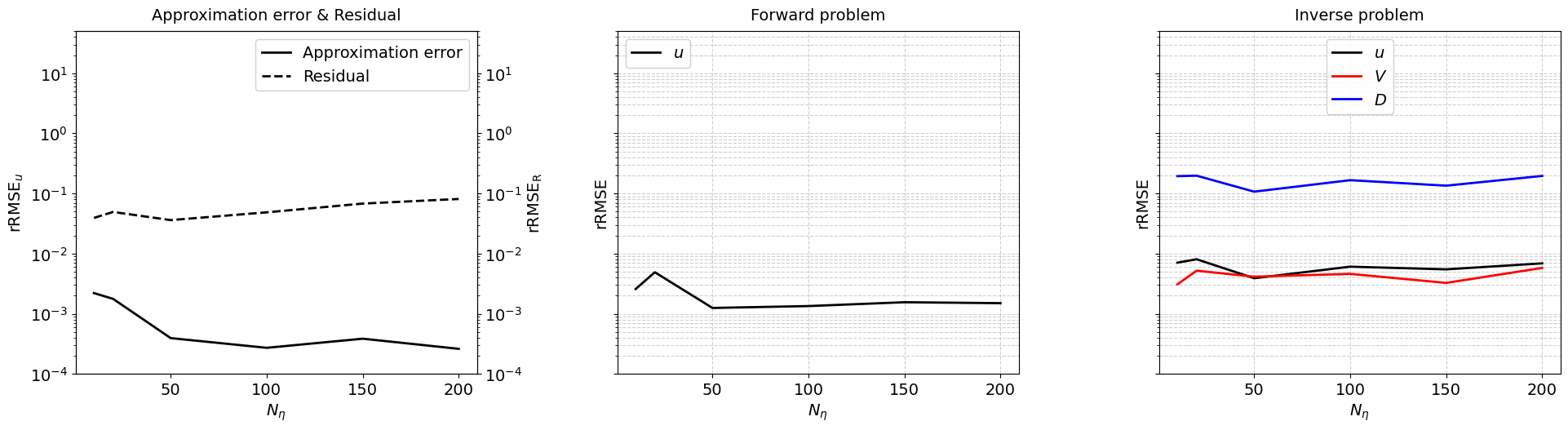}
    \caption{
ADE problem: 
\textbf{Left:} Approximation error of the pretrained fluctuation network and the corresponding residual error.
\textbf{Middle:} Relative error of the lPINN forward prediction. 
\textbf{Right:} Relative errors of the lPINN inverse prediction for the solution $u$ and the inferred parameters $V$ and $D$. 
$N_{\mathrm{train}}=500$.  The lPINN error for $N_{\eta}=50$, $N_{\mathrm{res}}=500$ in this test is not directly comparable with the lPINN error for $N_{\eta}=50$, $N_{\mathrm{res}}=500$ in the following forward problem as in that problem we use a shared $\lambda=10^{-4}$ across $N_{\mathrm{res}}$.
}
    \label{fig:Error_vs_num_eigenvalues_ADE}
\end{figure}

Figure~\ref{fig:Error_vs_num_eigenvalues_ADE} plots the relative approximation error, the residual error, and the relative errors in the forward and inverse lPINN solutions as functions of $N_\eta$, the number of neurons in the last hidden layers. The approximation error decreases as $N_\eta$ increases. The residual error stays nearly constant for $N_\eta<50$ and then increases with $N_\eta$. It should be noted that in the described experiment, the number of residual points in the offline step is fixed, which can cause the increase in the residual error with increasing $N_\eta$. 

The lPINN forward solution error increases with $N_\eta$ when the regularization coefficient $\lambda$ in Eq. \eqref{eq:RLS_final} is set to zero. For the considered linear ADE, Eq. \eqref{eq:RLS_final} yields a linear least-square problem that can be written as $\bm{A} \cdot \bm{W}_{N+1} = \bm{b}$, where $\bm{A}$ is a $N_{\mathrm{res}}\times {N_\eta}$ matrix and $\bm{b}$ is a $N_{\mathrm{res}}$ vector. 
The unregularized least-squares solution is unique when $\bm{A}$ has full column rank, equivalently when $\bm{A}^T \bm{A}$ is nonsingular. Regularization, i.e., setting $\lambda>0$, can provide uniqueness and improve numerical stability when the system is rank deficient or ill-conditioned.
We select $\lambda$ for each $N_\eta$ through a grid search. We do not observe a decrease in the lPINN rRMSE for $N_\eta>30$ because the residual rRMSE increases with $N_\eta$. We show in Section \ref{sec:superresolution} that the pretraining accuracy can be improved by increasing the size of the training dataset or, as stated earlier, using more residual points at the offline step.

For the inverse problem, the relative errors in the estimated $u$, $V$, and $D$ are weakly dependent on $N_\eta$. The overall best performance is achieved for $N_{\eta}=50$, the value we use in the remaining numerical experiments of this section.  

In the remainder of this section, we compare lPINN with vanilla PINN and PICKLE ROM.

\subsubsection{Comparison with PINN}
\begin{figure}[H]
    \centering
        \includegraphics[width=\linewidth]{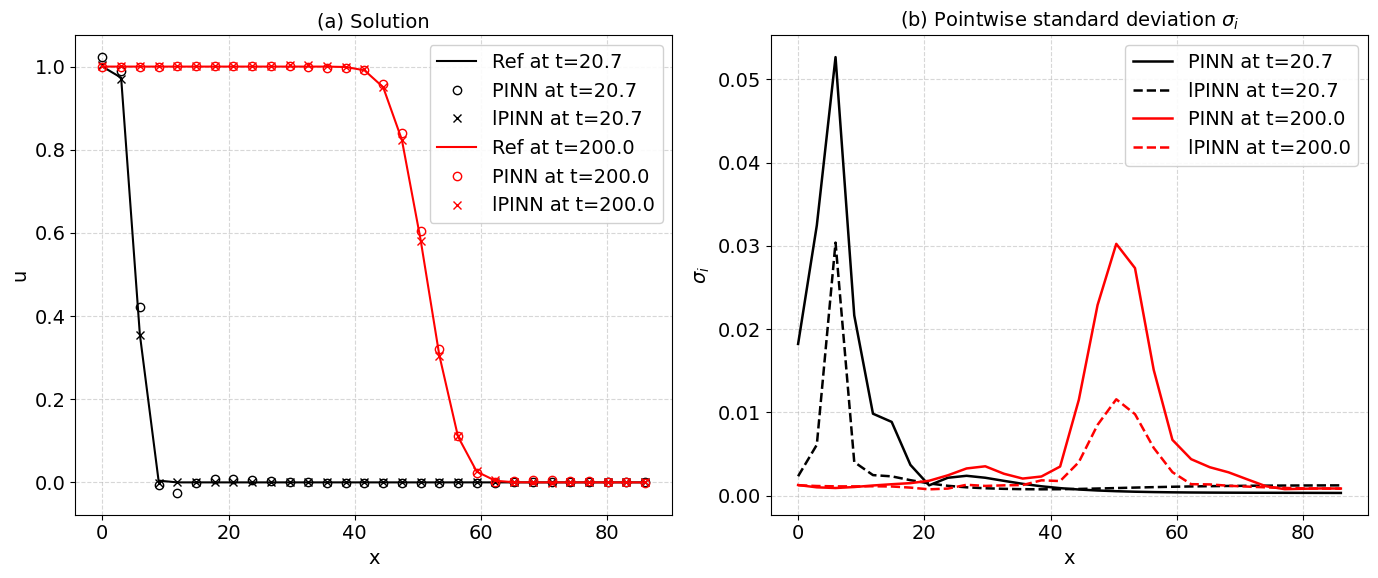}
    \caption{
    ADE forward problem: (a) Comparison of the PINN and lPINN predictions with the reference solution at $t=20.7$ s, and $200.0$ s. (b) Pointwise standard deviation at the $N_x=30$ grid points on the same mesh at $t=20.7$ s, and $200.0$ s. $N_{\mathrm{res}}=500$, $N_\eta=50$, $N_{\mathrm{train}}=500$. 
    }
    \label{fig:pinn_vs_lPINN_ADE_forward}
\end{figure}

We first consider the forward problem. Figure~\ref{fig:pinn_vs_lPINN_ADE_forward} compares the PINN and lPINN predictions against the reference solution at $t=20.7$ and $200.0$ s using $N_{\mathrm{res}}=500$ residual points. The lPINN solution is visibly closer to the reference and achieves an rRMSE of $0.0015$, whereas PINN yields an rRMSE of $0.0104$. The computational cost is also significantly reduced: lPINN requires $0.2$ s, while PINN requires $631.2$ s. All PINN and lPINN computations are performed on the same workstation using an NVIDIA RTX A6000 GPU with 48 GB of memory.

\begin{figure}[H]
    \centering
      \includegraphics[width=\linewidth]{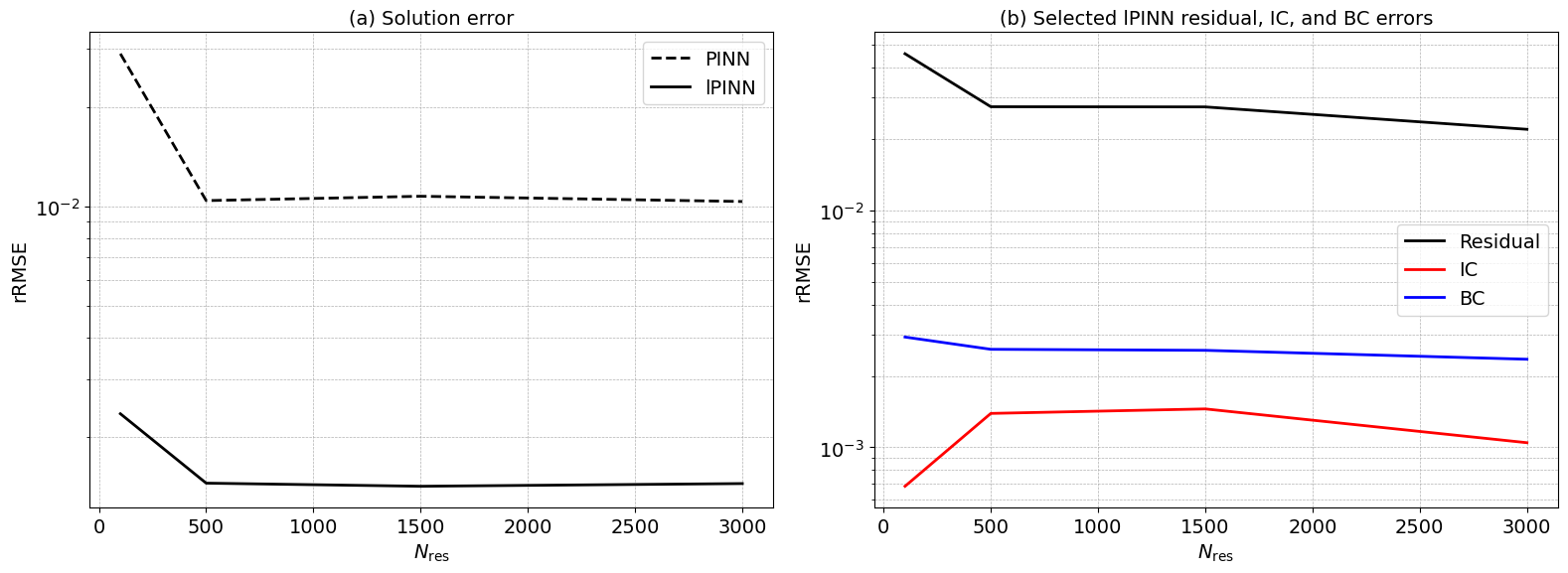}
   
\caption{
ADE forward problem: (a)  rRMSE$(\overline{\bm{u}})$ as a function of the number of residual points $N_{\mathrm{res}}$, where $\overline{\bm{u}}$ is the mean of 10 PINN or lPINN solutions. (b)  Relative root-mean-square error (rRMSE) of residual, initial condition, and boundary condition error as a function of the number of residual points $N_{\mathrm{res}}$. $N_\eta=50$, $N_{\mathrm{train}}=500$. 
}
    \label{fig:Error_forward_ADE}
\end{figure}

We next examine the sensitivity of lPINN and PINN errors to the number of residual points. Figure~\ref{fig:Error_forward_ADE}(a) shows that lPINN rRMSEs are consistently lower than PINN's. 
Both lPINN and PINN errors decrease with increasing $N_{\mathrm{res}}$,  and reach an asymptotic value at $N_{\mathrm{res}}=500$. 
 
\begin{figure}[H]
    \centering
        \includegraphics[width=\linewidth]{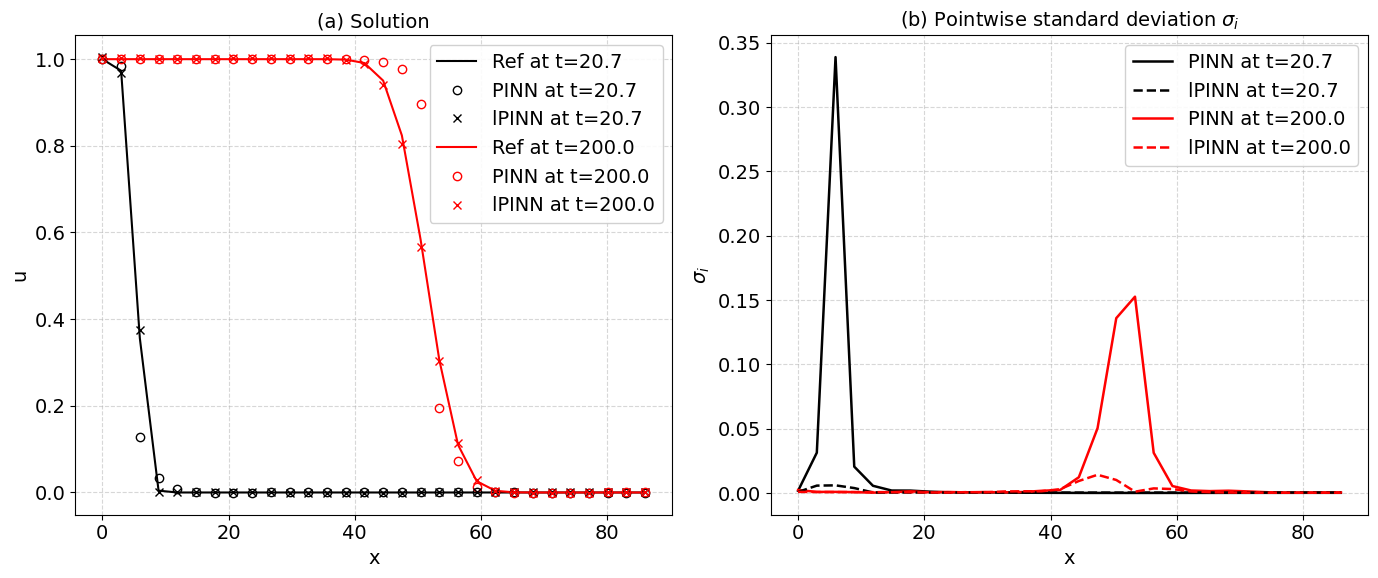}
\caption{
ADE inverse problem: Comparison of the PINN and lPINN estimations of $u(x,t)$  with the reference solution at $t=20.7$ s, and $200.0$ s.
$N_{\mathrm{res}}=500$, $N_{\mathrm{m}}=40$, $N_\eta=50$, $N_{\mathrm{train}}=500$.
}
    \label{fig:pinn_vs_lPINN_ADE_inv}
\end{figure}

Next, we compare PINN and lPINN inverse solutions for the parameters $V$ and $D$ and the state $u$ obtained with $N_{\mathrm{m}}=40$ measurements of $u$.  Figure~\ref{fig:pinn_vs_lPINN_ADE_inv} shows estimated $u$ and the reference $u$ at times $t=20.7$ and $200.0$ s. The lPINN prediction achieves an rRMSE of $0.0039$, compared with $0.0509$ for PINN. The $V$ and $D$ parameters found with lPINN have rRMSE values of $0.0041$ and  $0.1074$, respectively.  PINN yields $V$ with an rRMSE of $0.0290$ and $D$ with an rRMSE of $0.8513$. 
The inference time is reduced from $350.6$ s for PINN to $19.2$ s for lPINN.

\begin{figure}[H]
    \centering
    \includegraphics[width=0.75\linewidth]{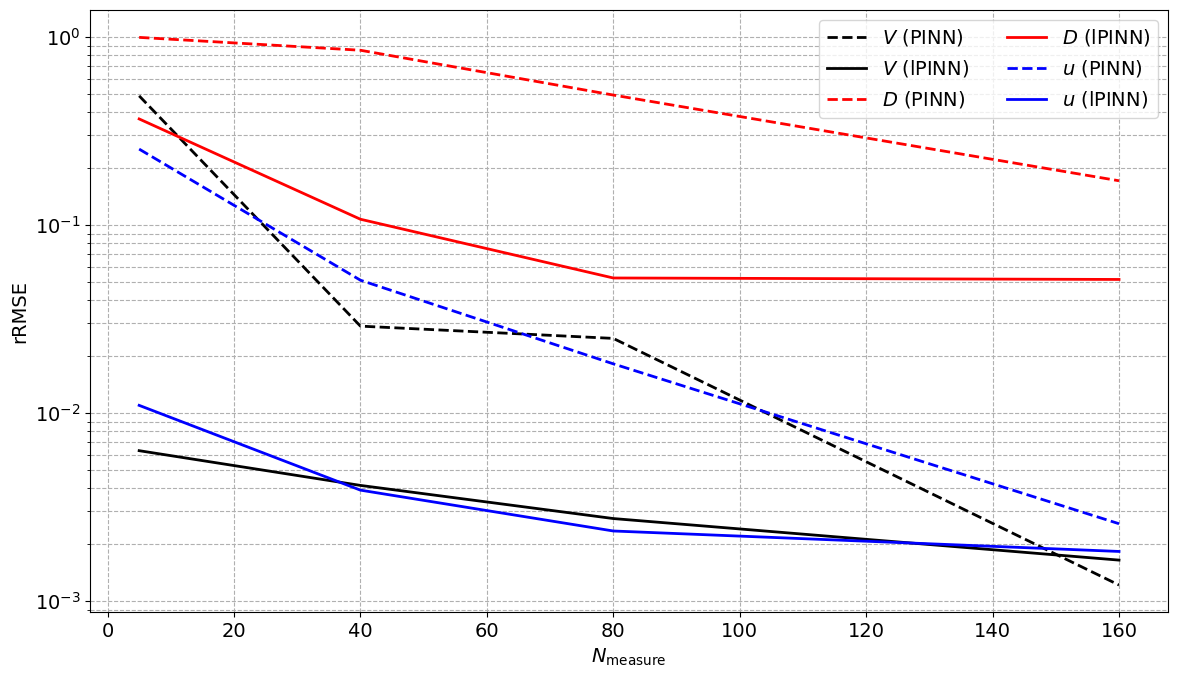}
\caption{
ADE inverse problem: Relative root-mean-square error (rRMSE) of the inferred velocity $V$, dispersion coefficient $D$, and predicted solution $u$ as functions of the number of measurement points $N_{\mathrm{m}}$.  $N_{\mathrm{res}}=500$, $N_\eta=50$, $N_{\mathrm{train}}=500$.
}
    \label{fig:Error_inv_ADE}
\end{figure}

Figure~\ref{fig:Error_inv_ADE} shows the rRMSEs in the lPINN and PINN inverse solutions as functions of $N_{\mathrm{m}}$. The errors are reported for the estimated $V$, $D$, and $u$. For all three variables, lPINN generally outperforms PINN, with the largest advantage occurring for smaller $N_{\mathrm{m}}$. This behavior matters in practical inverse problems, where observations are often limited and expensive to obtain.

\subsubsection{Comparison with the dPICKLE ROM method.}

Here, we compare the forward lPINN and dPICKLE solutions. The dPICKLE ROM method \cite{tartakovsky2024physics} is based on space--time-dependent Karhunen--Loève expansions (KLEs) of the state variable $u$:
\begin{equation}
\bm{u}(\bm{a})
=
\bar{\bm{u}}
+
\bm{\Phi}\bm{a},
\label{eq:dpickle_kle}
\end{equation}
where $\bm{u}$ is the solution vector, $\bar{\bm{u}}$ is the training dataset mean,  $\bm{\Phi}$ is the matrix of eigenvectors scaled with the square root of the corresponding eigenvalues, and $\bm{a}$ is the vector of coefficients. These coefficients can be estimated using projection, or, as in PINN and lPINN methods, by minimizing the L2 norm of the PDE residuals:
\begin{equation}
\bm{a}^{*}
=
\arg\min_{\bm{a}}
|| \bm{R} (\bm{a}) ||_2^2
\label{eq:dpickle_burgers}
\end{equation}
where $\bm{R} (\bm{a})$ is the vector of residuals. In dPICKLE, as in most ROM methods, all quantities, including $\bar{\bm{u}}$, $\bm{\Phi}$, $\bm{R}$, and, ultimately, the solution $\bm{u}$ are computed on the space-time mesh as the one used in the training dataset generation.  
In \cite{tartakovsky2024physics}, the spatial and temporal derivatives entering $\bm{R}$ are evaluated using a finite difference discretization. 

For a fair comparison, we use the same number of samples to compute $\bar{\bm{u}}$ and $\bm{\Phi}$ as in the lPINN solution. Also, the dimensionality of $\bm{a}$ is set the same as that of $\bm{W}_{N+1}$. The residual points in the lPINN solution are chosen to coincide with the grid points where dPICKLE residuals are evaluated.   

The forward solutions are compared for the parameter values
\[
V=1.2V^{*},
\qquad
D=0.7D^{*}.
\]

\begin{figure}[H]
    \centering
    \includegraphics[width=0.5\linewidth]{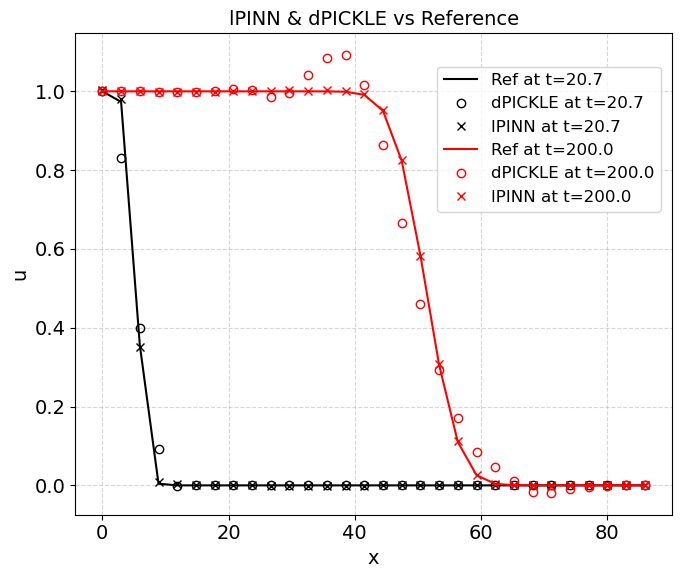}
    \caption{
    ADE forward problem: Comparison of the lPINN and dPICKLE predictions with the reference solution at  $t=20.7$ s and $200.0$ s. Here, $N_{\mathrm{res}}=784$, $N_\eta=50$, and $N_{\mathrm{train}}=500$.
    }
    \label{fig:lpinn_vs_dPICKLE_ADE_forward}
\end{figure}

Figure~\ref{fig:lpinn_vs_dPICKLE_ADE_forward} compares the lPINN and dPICKLE predictions with the reference solution at $t=20.7$ and $200.0$ s. The lPINN prediction is in closer agreement with the reference solution and achieves an rRMSE of $0.0013$, whereas dPICKLE yields an rRMSE of $0.0470$. However, dPICKLE has a lower computational cost: lPINN requires only $0.26$ s versus $0.01$ s for dPICKLE. The lPINN runtime is longer because it takes more time to compute derivatives via automatic differentiation than numerically in dPICKLE.  

\subsection{Burgers' equation}
\label{Burger's equation numerical example}

In this section, we consider the one-dimensional Burgers' equation,
\begin{equation}
\label{eq:Burgers_eq}
\frac{\partial u}{\partial t}
+
u\frac{\partial u}{\partial x}
=
\nu \frac{\partial^2 u}{\partial x^2},
\end{equation}
subject to the initial condition
\begin{equation}
u(x,0)=-\sin(\pi x),
\end{equation}
and the boundary conditions
\begin{equation}
u(-1,t)=u(1,t)=0.
\end{equation}

The forward lPINN model is pretrained to solve the Burgers' equation on the time domain $(0,1]$ for $\nu \in \left[\frac{0.1}{\pi},\frac{0.5}{\pi}\right]$. 
We construct the training dataset using the Cole-Hopf analytical solution of the Burgers' equation, 
\begin{equation}
u(x,t)=
\frac{
\displaystyle \int_{-\infty}^{\infty}
-\sin\!\bigl(\pi(x-\xi)\bigr)\,
\exp\!\left(-\frac{\cos(\pi(x-\xi))}{2\pi \nu}\right)
\exp\!\left(-\frac{\xi^2}{4\nu t}\right)\, d\xi
}{
\displaystyle \int_{-\infty}^{\infty}
\exp\!\left(-\frac{\cos(\pi(x-\xi))}{2\pi \nu}\right)
\exp\!\left(-\frac{\xi^2}{4\nu t}\right)\, d\xi
}.
\end{equation}
We uniformly draw $N_{\mathrm{train}}=500$ samples of $\nu$ from the specified interval.
For each $\nu$ sample, the Cole-Hopf solution is evaluated on an $N_x \times N_t$ mesh with $N_x=30$ and $N_t=30$. The temporal grid is uniform, and the spatial grid is center-clustered over $[-1,1]$. The center-clustered spatial grid points are defined as,
\[
s_i=-1+\frac{2i}{N_x-1},
\qquad
x_i=\operatorname{sign}(s_i)
\left[
1-\frac{\tanh\!\left(\sigma(1-|s_i|)\right)}
{\tanh(\sigma)}
\right],
\qquad i=0,\ldots,N_x-1.
\]
where clustering strength is set to $\sigma=2$.

The pretrained forward lPINN model is then tested for
\[
\nu^{\text{ref}}=\frac{0.1}{\pi},
\]
using the same mesh. The inverse problem is solved for $\nu$ and $u$ given sparse measurements of $u$.

The rRMSEs in the predicted solution $\mathbf{u}$ evaluated on the $N_x \times N_t$ mesh and, for the inverse problem, estimated viscosity $\nu$ are defined as
\[
\mathrm{rRMSE}_{\mathbf{u}}
=
\frac{
\left\|
\mathbf{u}-\mathbf{u}^{\mathrm{ref}}
\right\|_2
}{
u_0\sqrt{N_u}
},
\qquad
\mathrm{rRMSE}_{\nu}
=
\frac{
\left|\nu-\nu^{\mathrm{ref}}\right|
}{
\left|\nu^{\mathrm{ref}}\right|
}.
\]
where $N_u=N_x \times N_t$, $\bm{u}^{\text{ref}}$ is the reference solution evaluated on the $N_x \times N_t$ mesh, and $u_0 = u(x=-\frac{1}{2},t=0)$ is the maximum of the initial condition. The dimensionless rRMSEs for the PDE residual, initial condition, and boundary condition are defined as
\[
\mathrm{rRMSE}_R
=
\frac{L}{u_0^2}
\frac{\left\|\bm{\mathcal{R}}\right\|_2}{\sqrt{N_{\mathrm{res}}}},
\qquad
\mathrm{rRMSE}_{\mathrm{IC}}
=
\frac{
\left\|
\mathbf{u}_{\mathrm{IC}}-\mathbf{u}_{\mathrm{IC}}^{\mathrm{ref}}
\right\|_2
}{
u_0\sqrt{N_{\mathrm{IC}}}
},
\]
and
\[
\mathrm{rRMSE}_{\mathrm{BC}}
=
\frac{
\left\|
\mathbf{u}(x=-1)-\mathbf{u}^{\mathrm{ref}}(x=-1)
\right\|_2
}{
u_0\sqrt{N_{\mathrm{BC},-1}}
}
+
\frac{
\left\|
\mathbf{u}(x=1)-\mathbf{u}^{\mathrm{ref}}(x=1)
\right\|_2
}{
u_0\sqrt{N_{\mathrm{BC},1}}
}.
\]

\begin{figure}[H]
        \includegraphics[width=1\linewidth]{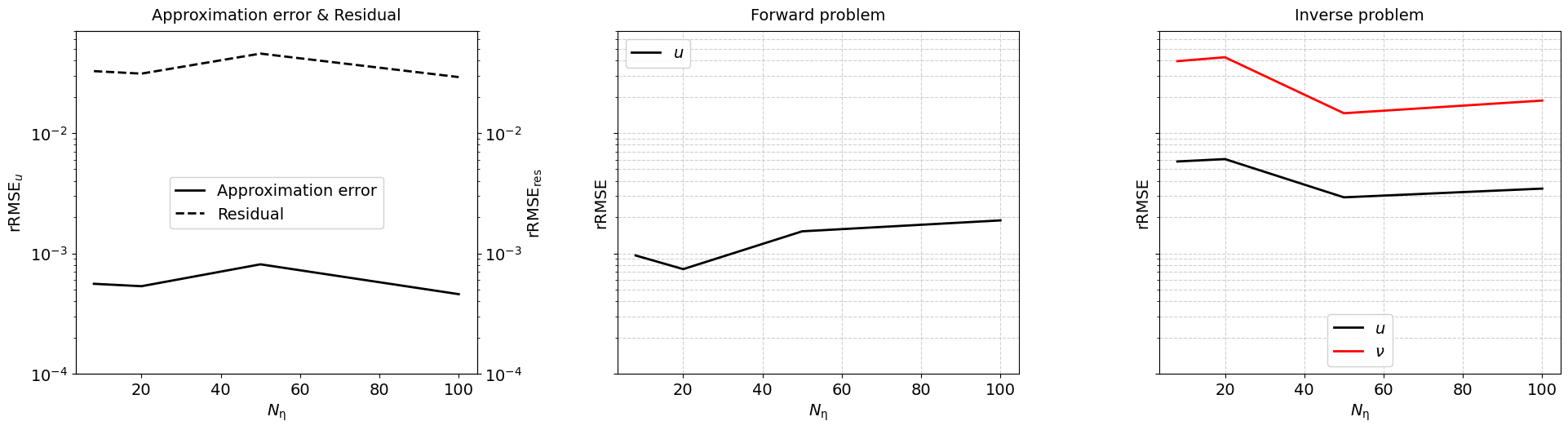}
    \caption{
    Burgers' equation: 
\textbf{Left:} Approximation error of the pretrained fluctuation network and the corresponding residual error.
\textbf{Middle:} Relative error of the lPINN forward prediction. 
\textbf{Right:} Relative errors of the lPINN inverse prediction for the solution $u$ and the inferred parameter  $\nu$.
Here, $N_{\mathrm{train}}=500$.
    }
    \label{fig:Error_vs_num_eigenvalues_burgers}
\end{figure}

Figure~\ref{fig:Error_vs_num_eigenvalues_burgers} plots the relative approximation error,  the residual of the approximated solution, and the relative errors in the forward and inverse lPINN solutions as functions of  $N_\eta$.  The number of residual points is set to  $N_{\mathrm{res}}=10N_\eta$. The inverse solution is obtained using  $N_{\mathrm{m}}=50$ measurements. The DNNs are pretrained using residual-based loss functions. This choice is justified by the comparison study summarized in Table \ref{tab:burgers_pretrain_compare}.
The approximation errors of the pretrained DNNs and the residual RMSE are weakly dependent on the considered $N_\eta$ values. The smallest rRMSE of the forward lPINN solution is obtained for $N_\eta=20$. In the inverse problem, the errors in both the estimated $u$ and $\nu$ are smallest at $N_\eta=50$. 
Accordingly, in the remaining examples of this section, we set $N_{\eta}=20$ for forward problems and $N_{\eta}=50$ for the inverse problems.

Next, we compare derivative-matching pretraining with residual-based pretraining. 
For residual-based offline pretraining, we split the $N_x \times N_t=900$ solution-grid points into 10 batches at each epoch. For each batch, the solution values of all $N_{\mathrm{train}}=500$ realizations at 90 grid points are used for training, together with $N_{\mathrm{res}}=90$ residual points, $N_{\mathrm{IC}}=90$ initial condition collocation points, and
$N_{\mathrm{BC}}=90$ boundary condition collocation points are randomly resampled for each batch.
For derivative-matching offline pretraining, the same 900 solution-grid points are split into 10 batches, and for each batch, the solution and derivative values at the corresponding 90 grid points are used for training, together with $N_{\mathrm{IC}}=90$ and $N_{\mathrm{BC}}=90$ randomly resampled initial condition and boundary condition collocation points.
The number of residual points is set to $N_{\mathrm{res}}=100$ in both lPINN and PINN for the forward problem, and we use $N_{\mathrm{res}}=100$ residual points and $N_{\mathrm{m}}=5$ measurements of $u$ for the inverse problem. 
As shown in Table~\ref{tab:burgers_pretrain_compare}, residual-based pretraining gives smaller approximation errors for both $N_\eta=20$ and $N_\eta=50$. Residual-based pretraining gives smaller residual RMSE and better accuracy in both the forward and inverse problems. Specifically, it reduces the forward error from $3.15\times10^{-3}$ to $8.07\times10^{-4}$ and also gives lower inverse errors for both $u$ and $\nu$. In the remaining numerical experiments of this section, we use the residual-based pretraining.

\begin{table}[H]
\centering
\caption{Comparison of residual-based and derivative-matching pretraining for Burgers' equation. As we use $N_{\mathrm{m}}=5$ in this test and we use $N_{\mathrm{m}}=50$ for solving the inverse problem in the following part of the section, the result is not directly comparable.  }
\label{tab:burgers_pretrain_compare}
\begin{tabular}{lccc}
\hline
Metric & $N_\eta$ & Residual-based & Derivative-matching \\
\hline
Approximation rRMSE
    & 20 & $\mathbf{5.35\times10^{-4}}$ & $2.11\times10^{-3}$ \\
    & 50 & $\mathbf{8.13\times10^{-4}}$ & $8.24\times10^{-4}$ \\
\hline
Residual rRMSE& 20 & $\mathbf{3.10\times10^{-2}}$ & $4.90\times10^{-2}$ \\
    & 50 & $\mathbf{4.56\times10^{-2}}$ & $6.32\times10^{-2}$ \\
\hline
Forward rRMSE$_u$
    & 20& $\mathbf{8.07\times10^{-4}}$& $3.15\times10^{-3}$\\
Inverse rRMSE$_u$
    & 50& $\mathbf{1.57\times10^{-3}}$ & $3.24\times10^{-3}$ \\
Inverse rRMSE$_\nu$
    & 50& $\mathbf{3.18\times10^{-3}}$ & $1.14\times10^{-2}$ \\
\hline
\end{tabular}
\end{table}

\begin{figure}[H]
    \centering
     \includegraphics[width=\linewidth]{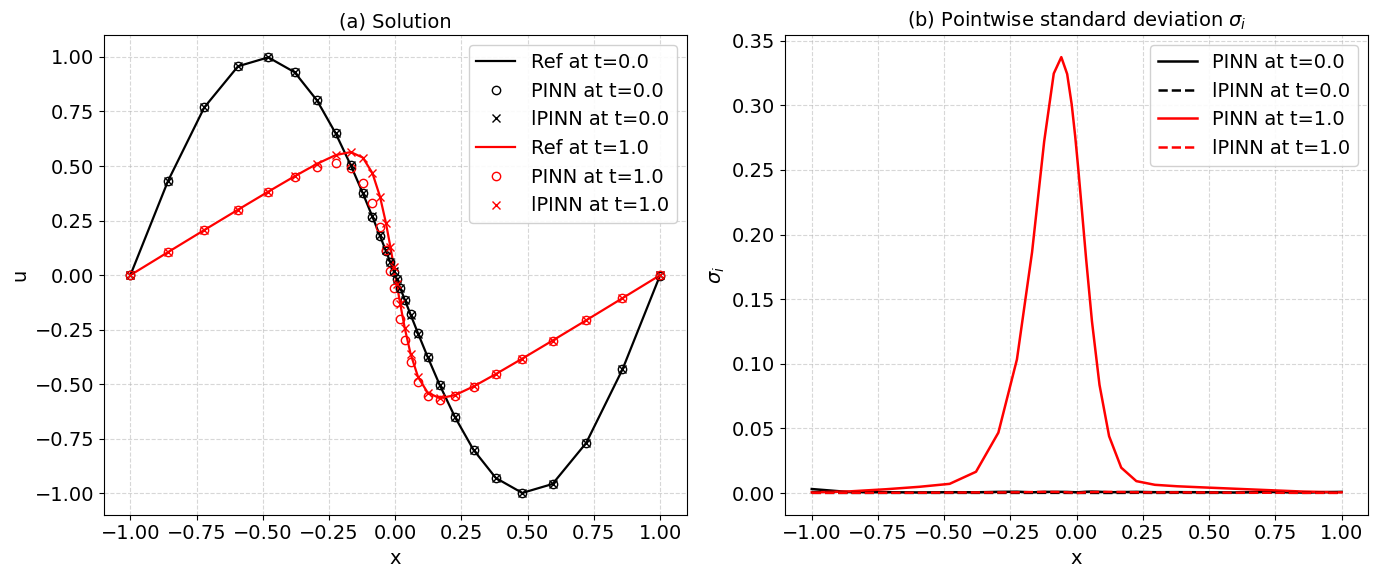}
    \caption{
    Burgers' forward problem: Comparison of the PINN and lPINN predictions with the reference solution at $t=0.0$ s and $1.0$ s. Here, $N_{\mathrm{res}}=100$, $N_\eta=20$, $N_{\mathrm{train}}=500$
    }
    \label{fig:pinn_vs_lPINN_burgers_forward}
\end{figure}

Next, we compare the PINN and lPINN solutions of the forward Burgers' problem. Figure~\ref{fig:pinn_vs_lPINN_burgers_forward} shows the PINN, lPINN, and reference solutions at $t=0.0$ and $1.0$ s. 
The lPINN solution agrees much more closely with the reference and achieves an rRMSE of $0.0008$, whereas the PINN yields an rRMSE of $0.0439$. The computational advantage is also substantial: lPINN requires only $17.3$ s, while PINN requires $308.3$ s.

\begin{figure}[H]
    \centering
       \includegraphics[width=1\linewidth]{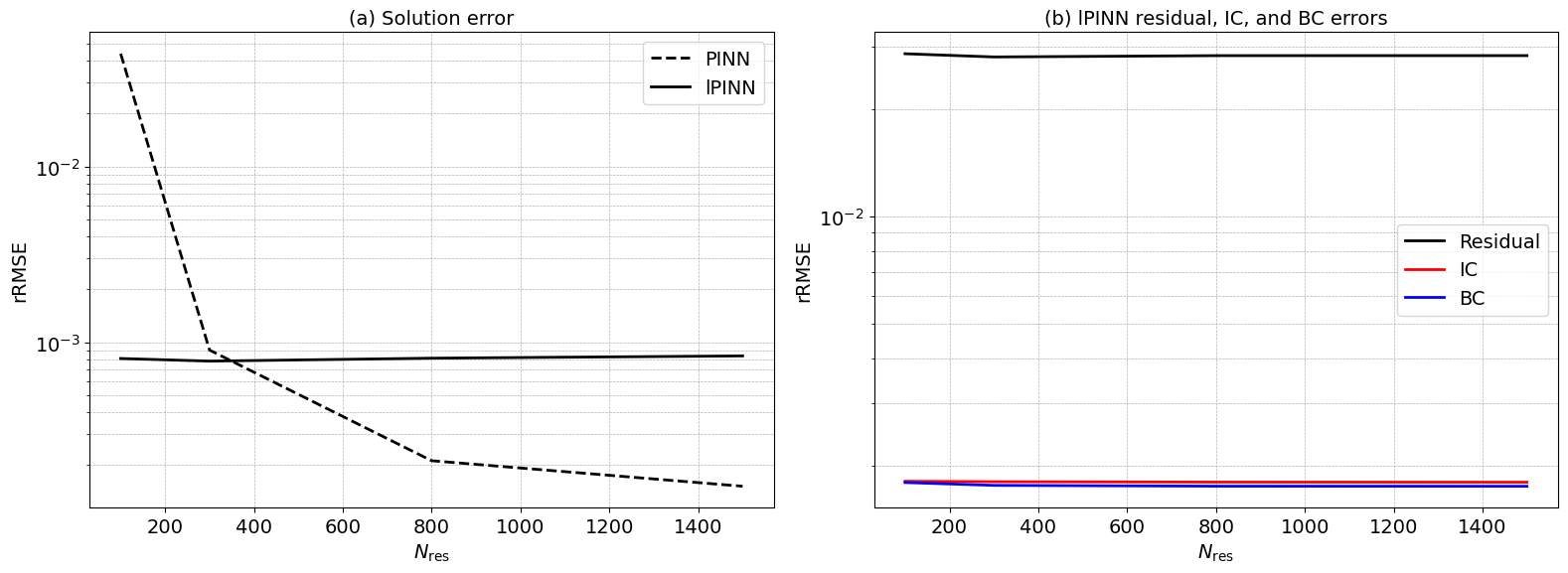}
    \caption{
Burgers' forward problem: (a) Relative root-mean-square error (rRMSE) of the predicted solution for PINN and lPINN as a function of the number of residual points $N_{\mathrm{res}}$. (b)  Relative root-mean-square error (rRMSE) of the residual, initial condition, and boundary condition errors as functions of the number of residual points $N_{\mathrm{res}}$. $N_\eta=20$, $N_{\mathrm{train}}=500$. }
    \label{fig:Error_forward_burgers}
\end{figure}

Figure~\ref{fig:Error_forward_burgers}(a) examines the dependence of lPINN and PINN errors in the forward Burgers' problem on $N_{\mathrm{res}}$ for $N_\eta=20$. 
The rRMSE of PINN decreases to an asymptotic value at $N_{\mathrm{res}}\approx 800$, while the rRMSE of lPINN stays nearly constant across $N_{\mathrm{res}}$.
The lPINN asymptotic error is close to the approximation error for $N_{\eta}=20$, which is the lower limit of lPINN error. The PINN asymptotic error is smaller than that of lPINN. This is different from the ADE problem where lPINN remained more accurate than PINN for all considered $N_{\mathrm{res}}$. 
Figure~\ref{fig:Error_forward_burgers}(b) explains why the lPINN error does not decrease with increasing $N_{\mathrm{res}}$: while the PDE residual error decreases,  the initial and boundary condition residual errors stay unchanged with increasing $N_{\mathrm{res}}$.  

\begin{figure}[H]
    \centering
        \includegraphics[width=\linewidth]{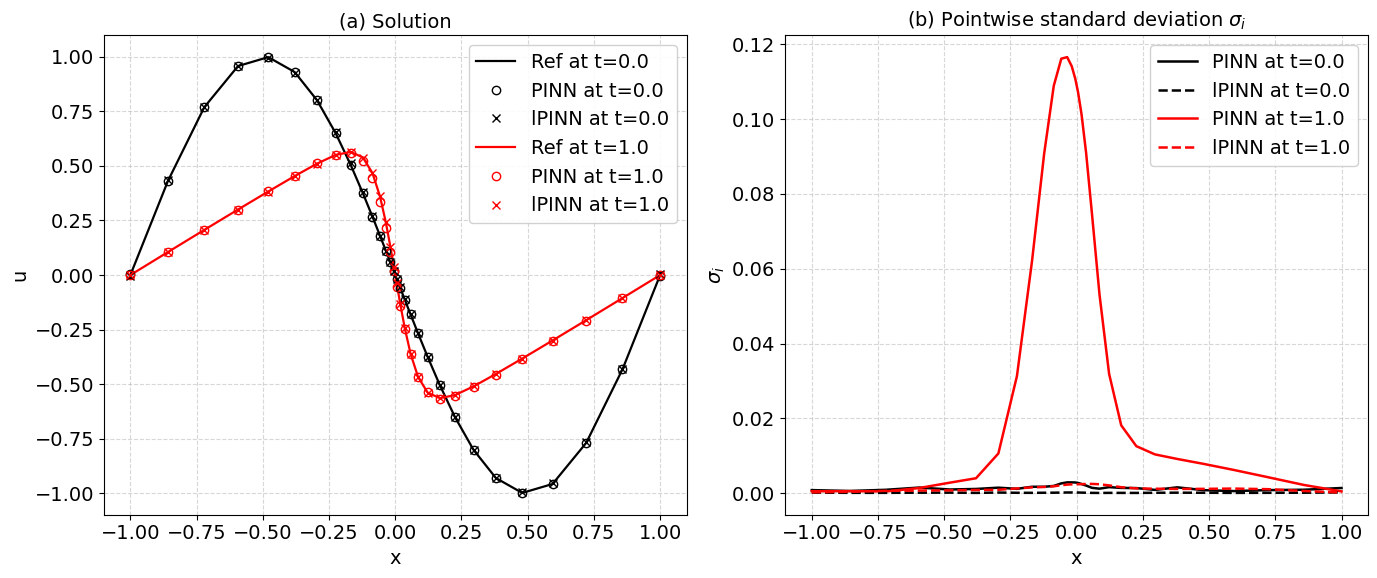}
    \caption{
    Burgers' inverse problem: Comparison of the PINN and lPINN predictions with the reference solution at $t=0.0$ and $1.0$ s. $N_{\mathrm{res}}=100$, $N_{\mathrm{m}}=5$, $N_\eta=50$, $N_{\mathrm{train}}=500$
    }
    \label{fig:pinn_vs_lPINN_burgers_inv}
\end{figure}

Figure~\ref{fig:pinn_vs_lPINN_burgers_inv} compares the PINN and lPINN $u$ inverse solutions at $t=0.0$, and $1.0$ s using $N_{\mathrm{res}}=100$ residual points and $N_{\mathrm{m}}=5$ measurements of $u$. The lPINN solution achieves an rRMSE of $0.0016$ versus $0.0061$ in PINN. The estimated $\nu$ in the lPINN solution has a relative error of $0.0032$ versus $0.0119$ in PINN. The time to obtain the lPINN solution (excluding pretraining) is $3.2$ s versus $251.5$ s in the PINN method. 

\begin{figure}[H]
    \centering
    \includegraphics[width=0.75\linewidth]{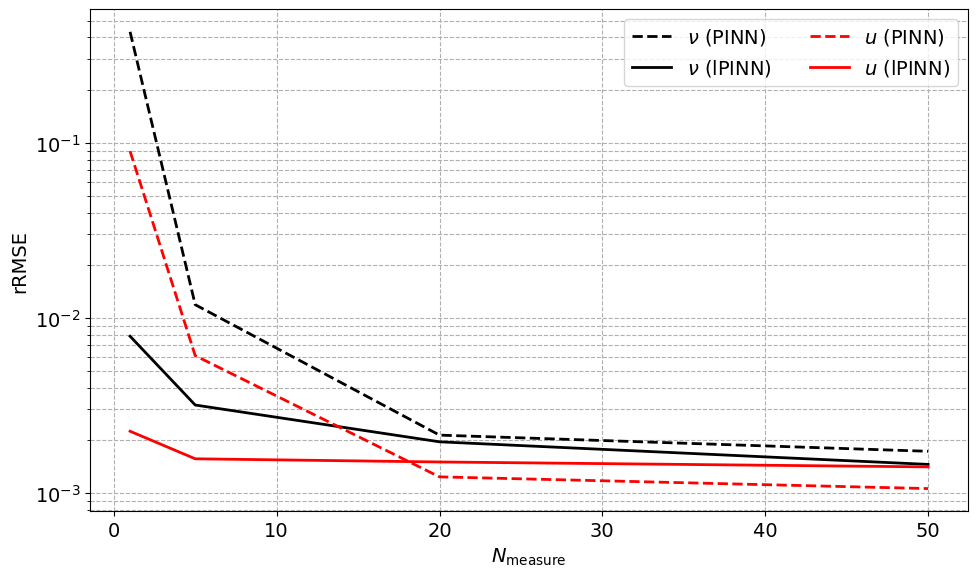}
    \caption{
    Burgers' inverse problem: Relative root-mean-square error (rRMSE) of the inferred kinematic viscosity $\nu$ and predicted solution $u$ as functions of the number of measurement points $N_{\mathrm{m}}$. $N_{\mathrm{res}}=100$, $N_\eta=50$, $N_{\mathrm{train}}=500$
    }
    \label{fig:Error_vs_Measurement_burgers}
\end{figure}

We further investigate the effect of the number of measurements on inverse solution accuracy for $N_\eta=50$ and $N_{\mathrm{res}}=100$. As shown in Figure~\ref{fig:Error_vs_Measurement_burgers}, the PINN errors in the estimated $\nu$ and  $u$ strongly depend on $N_{\mathrm{m}}$. In contrast, lPINN errors are less sensitive to $N_{\mathrm{m}}$. The lPINN method achieves similar accuracy with relatively few measurements. The error in the estimated $\nu$ is smaller than that of PINN for all considered $N_{\mathrm{m}}$. The $u$ estimate is more accurate in lPINN for $N_{\mathrm{m}} < 20$ and less accurate otherwise. This suggests that lPINN is more robust than PINN in the limited-data regime.

\subsection{Nonlinear Pendulum Equation}\label{sec:pendulum}

We next assess the proposed lPINN framework on the nonlinear pendulum equation
\begin{equation}
\label{eq:nlpend_eq}
\frac{d^2 \theta}{dt^2}
+
\gamma \frac{d\theta}{dt}
+
\frac{g_0}{\ell}\sin(\theta)
=0,
\end{equation}
subject to the initial conditions
\begin{equation}
\theta(0)=\frac{\pi}{2},
\qquad
\frac{d\theta}{dt}(0)=0,
\end{equation}
where $\gamma$ is the damping coefficient, $\ell$ is the pendulum length, and $\theta$ is the angle, $g_0$ is the gravity. 
In the forward problem, we compute $\theta$ for $\gamma\in[0.05,0.5]$, $\ell\in[0.5,2.0]$. In the inverse problem, $\gamma$, $\ell$, and $\theta(t)$ are estimated from $\theta$ measurements at a few time instances. 

To construct the training dataset, we draw $N_{\mathrm{train}}=500$ samples of $(\gamma,\ell)$ uniformly from their considered ranges.
 Eq.~\eqref{eq:nlpend_eq} is solved numerically for each parameter combination with the Runge–Kutta method on a uniform temporal grid using a constant time step $\Delta t = T/(N_t-1)$, where $N_t=300$ and $T=30$. In addition to $\theta$, $\frac{d\theta}{dt}$ and $\frac{d^2 \theta}{dt^2}$ are also obtained from the solver and used for the DNN pretraining. The pretrained model is then tested on the same temporal mesh for $\gamma^{\text{ref}}=0.1$ and $\ell^{\text{ref}}=0.8$.

Eq \eqref{eq:nlpend_eq} is challenging to solve using PINN because of $\theta$ oscillation with $t$ as result of the spectral bias, the tendency of DNNs to learn low-frequency components more readily than high-frequency components \cite{wang2021eigenvector}. In lPINN, we find that a standard fully connected network approximates $\theta(t)$ more accurately at early times than at later times. We also find that the training losses of the mean and fluctuation components converge at substantially different rates, which further complicates joint training with standard DNNs. 
To improve DNN pretraining, we separately pretrain the mean and fluctuation networks using the numerical values of $\bar{\theta}$, $d\bar{\theta}/dt$, $d^2\bar{\theta}/dt^2$, and $\theta'$, $d\theta'/dt$, $d^2\theta'/dt^2$, respectively. Also, we use a Fourier Feature Neural Network (FFNN) rather than a standard fully connected network to model fluctuations $\theta'(t)$. In both PINN and lPINN, we use FFNNs with a mapping size of 64 and the hidden-layer width $N_\eta$. 
We use different frequency scales (0.5 for PINN and 0.1 for lPINN) because preliminary tuning showed that the two methods achieve their best performance at different Fourier-feature scales.
For lPINN, the smaller scale provides sufficient resolution of the oscillatory solution without introducing spurious high-frequency artifacts.

The rRMSEs of the predicted $\bm{\theta}$ and, for the inverse problem, estimated $\gamma$ and $\ell$ are defined as
\[
\mathrm{rRMSE}_{\bm{\theta}}
=
\frac{
\left\|
\bm{\theta}-\bm{\theta}^{\mathrm{ref}}
\right\|_2
}{
\theta_0\sqrt{N_t}
},
\qquad
\mathrm{rRMSE}_{\gamma}
=
\frac{
\left|\gamma-\gamma^{\mathrm{ref}}\right|
}{
\left|\gamma^{\mathrm{ref}}\right|
},
\qquad
\mathrm{rRMSE}_{\ell}
=
\frac{
\left|\ell-\ell^{\mathrm{ref}}\right|
}{
\left|\ell^{\mathrm{ref}}\right|
},
\]
where $\bm{\theta}^{\mathrm{ref}}$ is the reference solution and  $\theta_0=\pi/2$ is the initial condition. The dimensionless rRMSEs for the PDE residual and the initial conditions are defined as
\[
\mathrm{rRMSE}_R
=
\frac{\ell}{\theta_0 g_0}
\frac{\left\|\bm{\mathcal{R}}\right\|_2}{\sqrt{N_{\mathrm{res}}}},
\qquad
\mathrm{rRMSE}_{\text{IC},0}
=
\frac{
\left|\theta(0)-\theta_0\right|
}{
\theta_0
},
\]
and
\[
\mathrm{rRMSE}_{\text{IC},1}
=
\frac{
\left|\frac{d{\theta}(0)}{dt}\right|
}{
\theta_0\sqrt{g_0/\ell}
}.
\]

For the nonlinear pendulum Eq \eqref{eq:nlpend_eq}, residual-based pretraining did not produce an accurate pretrained representation and resulted in large errors during the subsequent online step. In contrast, derivative-matching pretraining accurately reproduced the mean and fluctuation fields, along with their first and second time derivatives. Therefore, this pretraining approach is better suited for learning the basis functions for the nonlinear pendulum equation, and we use it in all examples in this section. 
For derivative-matching offline pretraining, the solution and its first- and second-order time derivatives of all $N_{\mathrm{train}}=500$ realizations at $N_t=300$ temporal grid points are used for training. 

\begin{figure}[H]
    \centering
        \includegraphics[width=1\linewidth]{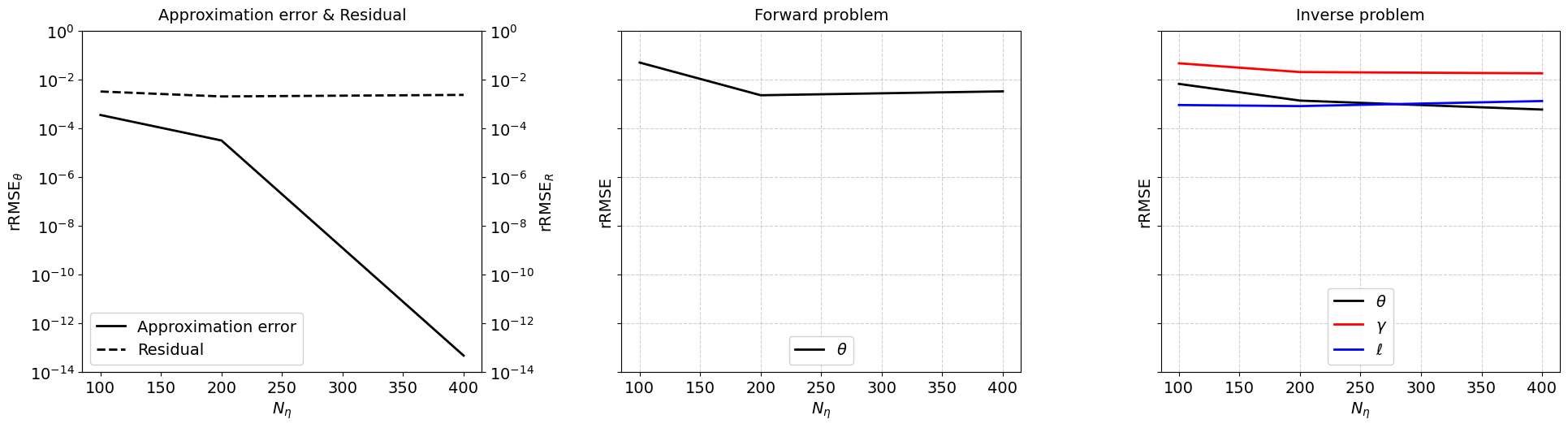}
    \caption{
    Nonlinear pendulum equation. 
\textbf{Left:} Approximation error of the pretrained fluctuation network and the corresponding residual error.
\textbf{Middle:} Relative error of the lPINN forward prediction. 
\textbf{Right:} Relative errors of the lPINN inverse prediction for the solution $\theta$ and the inferred parameter  $\gamma$ and $\ell$. 
Here, $N_{\mathrm{train}}=500$, $N_{\mathrm{res}}=10N_{\eta}$, and $N_{\mathrm{m}}=50$ for the inverse problem.}
    \label{fig:Error_vs_num_eigenvalues_nlpend}
\end{figure}

Figure~\ref{fig:Error_vs_num_eigenvalues_nlpend} plots the relative approximation error, the residual of the approximated solution, and the relative errors in the forward and inverse lPINN solutions as functions of $N_\eta$. The number of residual points is set to $N_{\mathrm{res}}=10N_\eta$ while the training set size is kept fixed at $N_{\mathrm{train}}=500$. For the inverse problem, the number of $\theta$ measurements is set to $N_{\mathrm{m}}=50$. The approximation error decreases significantly as $N_\eta$ increases, reaching nearly machine precision at $N_\eta=400$. In contrast, the residual error remains at approximately the same order of magnitude for all tested widths, indicating that increasing the last layer beyond $N_\eta=200$ does not improve the overall solution representation for the fixed $N_{\mathrm{train}}$. 
The forward lPINN solutions reach the smallest error at $N_{\eta}=200$. 
For the inverse problem, the relative errors in the estimated $\theta$ and the estimated parameter $\gamma$ decrease with $N_\eta$, while the estimated parameter $\ell$ is practically independent of it. 
Therefore, considering both accuracy and computational cost, we set $N_{\eta}=200$ in the following numerical experiments.

\begin{figure}[H]
    \centering
        \includegraphics[width=\linewidth]{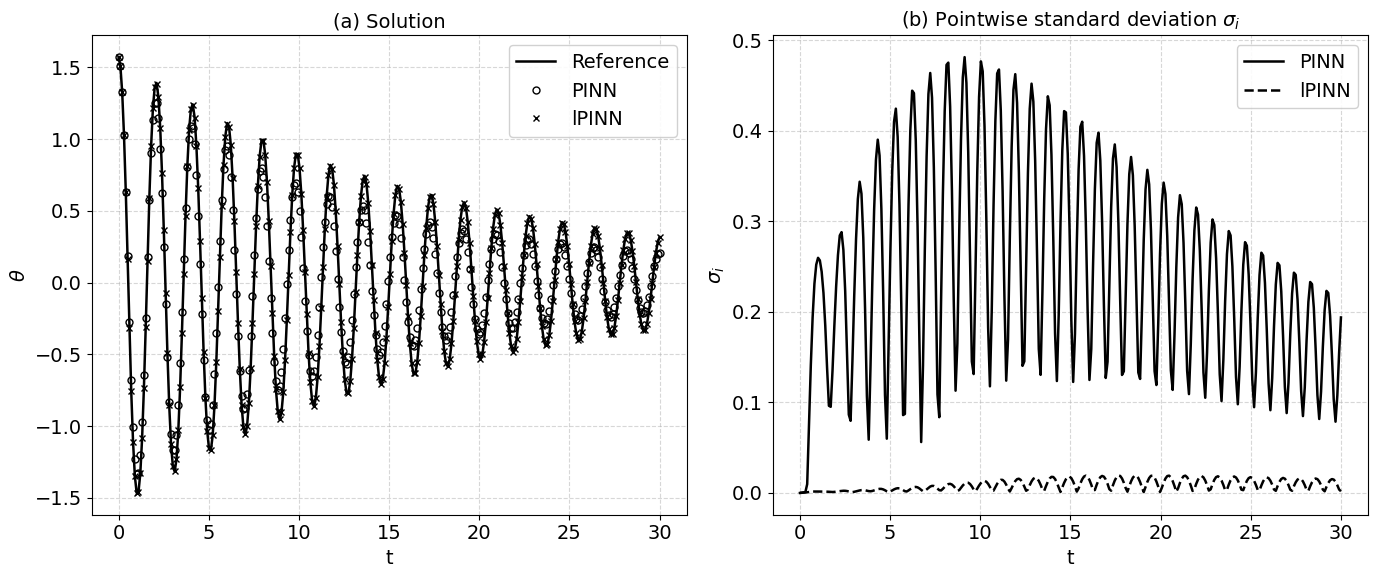}
    \caption{
    Nonlinear pendulum equation forward problem: Comparison of the PINN and lPINN predictions with the reference solution. $N_{\mathrm{res}}=1000$, $N_{\eta}=200$, $N_{\mathrm{train}}=500$.
    }
    \label{fig:pinn_vs_lPINN_nlpend_forward}
\end{figure}

Figure~\ref{fig:pinn_vs_lPINN_nlpend_forward} compares the PINN and lPINN forward solutions with the reference solution.  The lPINN achieves an rRMSE of $0.0022$, compared with $0.0994$ for PINN, corresponding to an error reduction of approximately $98\%$. The computational cost drops from $2375.3$~s for PINN to $287.1$~s for lPINN, an approximately $8.3$-fold speedup.

\begin{figure}[H]
    \centering
        \includegraphics[width=1\linewidth]{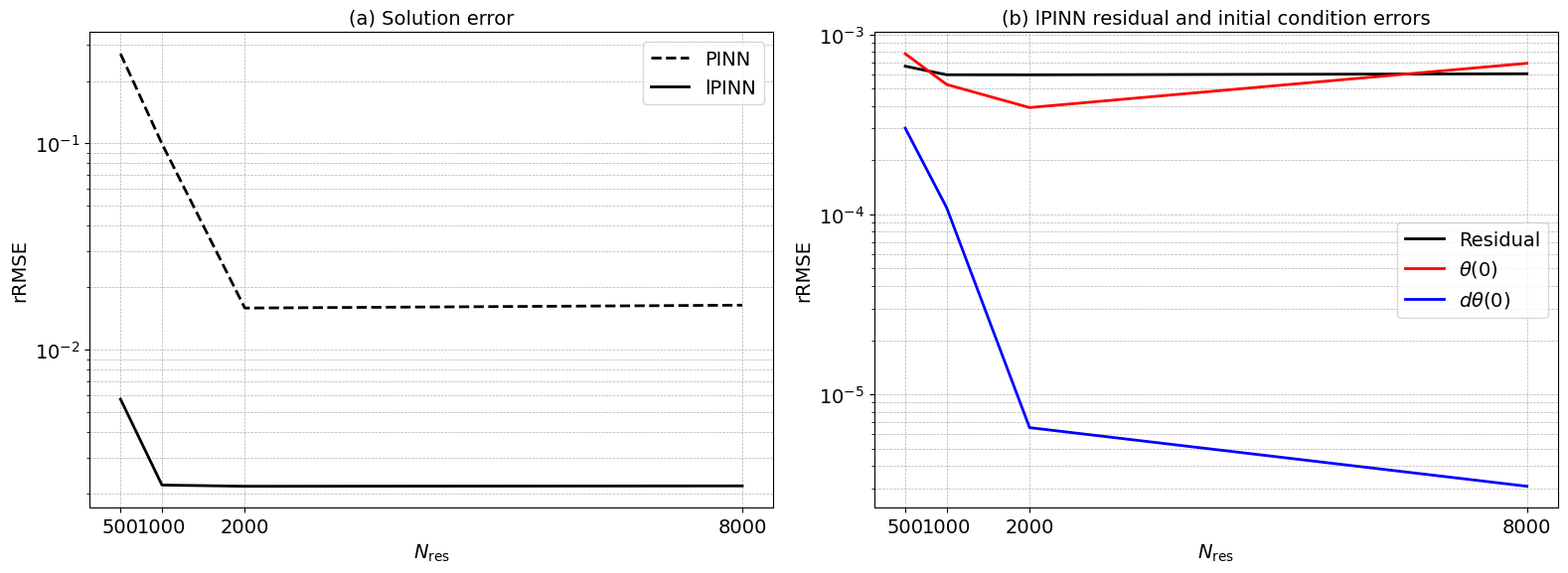}
    \caption{
    Nonlinear pendulum forward problem: (a) rRMSE of the predicted solution for PINN and lPINN as a function of the number of residual points $N_{\mathrm{res}}$. (b) rRMSE of residual and initial condition errors as a function of the number of residual points $N_{\mathrm{res}}$. 
    $N_{\eta}=200$, $N_{\mathrm{train}}=500$.
   }
    \label{fig:Error_forward_nlpend}
\end{figure}

We next examine how the number of residual points affects forward solution accuracy.  Figure~\ref{fig:Error_forward_nlpend}(a) shows that lPINN errors are smaller than PINN errors for all considered $N_{\mathrm{res}}$. The PINN error initially decreases with $N_{\mathrm{res}}$ and reaches an asymptotic value at $N_{\mathrm{res}}=2000$. 
The lPINN error also initially decreases with $N_{\mathrm{res}}$ and reaches an asymptotic value at $N_{\mathrm{res}}=1000$. 
At $N_{\mathrm{res}}=1000$, the lPINN error is more than an order of magnitude smaller than the PINN error.  

\begin{figure}[H]
    \centering
    \includegraphics[width=\linewidth]{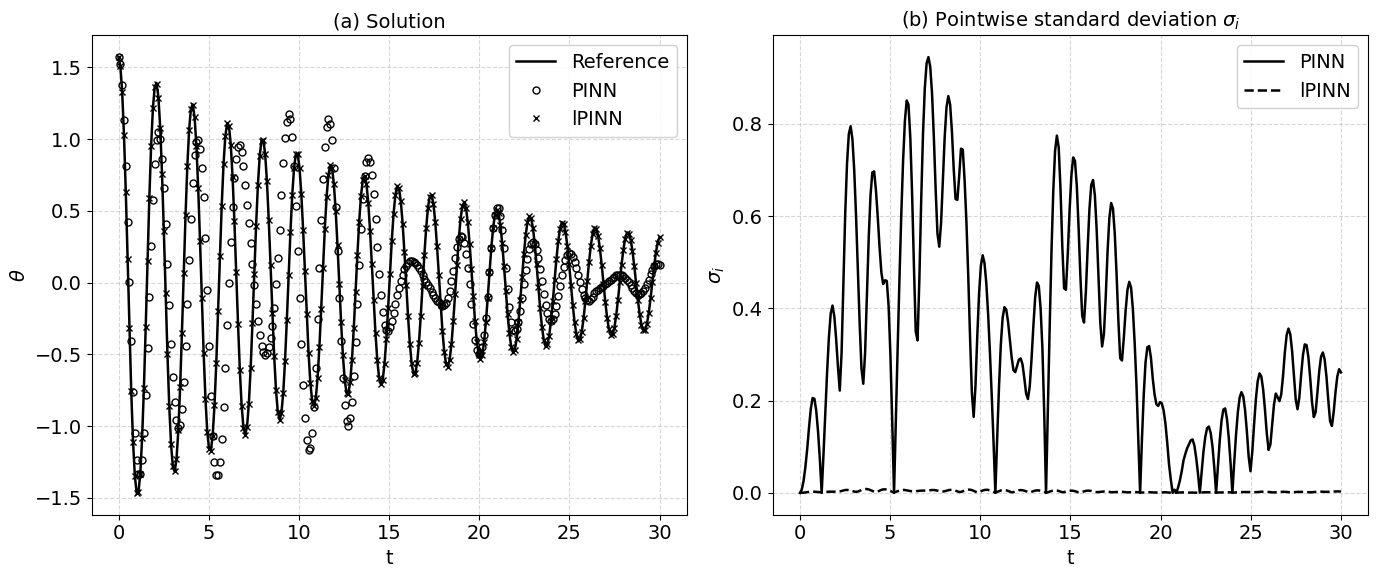}
    \caption{
    Nonlinear pendulum equation inverse problem: Comparison of the PINN and lPINN predictions with the reference solution. $N_{\mathrm{res}}=500$, $N_{\mathrm{m}}=10$, $N_{\eta}=200$, $N_{\mathrm{train}}=500$. 
    }
    \label{fig:pinn_vs_lPINN_nlpend_inv}
\end{figure}

Figure~\ref{fig:pinn_vs_lPINN_nlpend_inv} compares the inverse PINN and lPINN solutions for $\theta$ with the reference solution. The lPINN solution is more accurate, with an rRMSE of $0.0014$ versus $0.3767$ in the PINN solution. The lPINN parameter estimates are also more accurate, yielding rRMSEs of $0.0133$ for $\gamma$ and $0.0012$ for $\ell$ versus $1.0917$ and $0.2648$, respectively, in the PINN method. The online part of the lPINN solution takes $209.0$ s versus $2271.0$ s for the PINN solution.

\begin{figure}[H]
    \centering
    \includegraphics[width=0.75\linewidth]{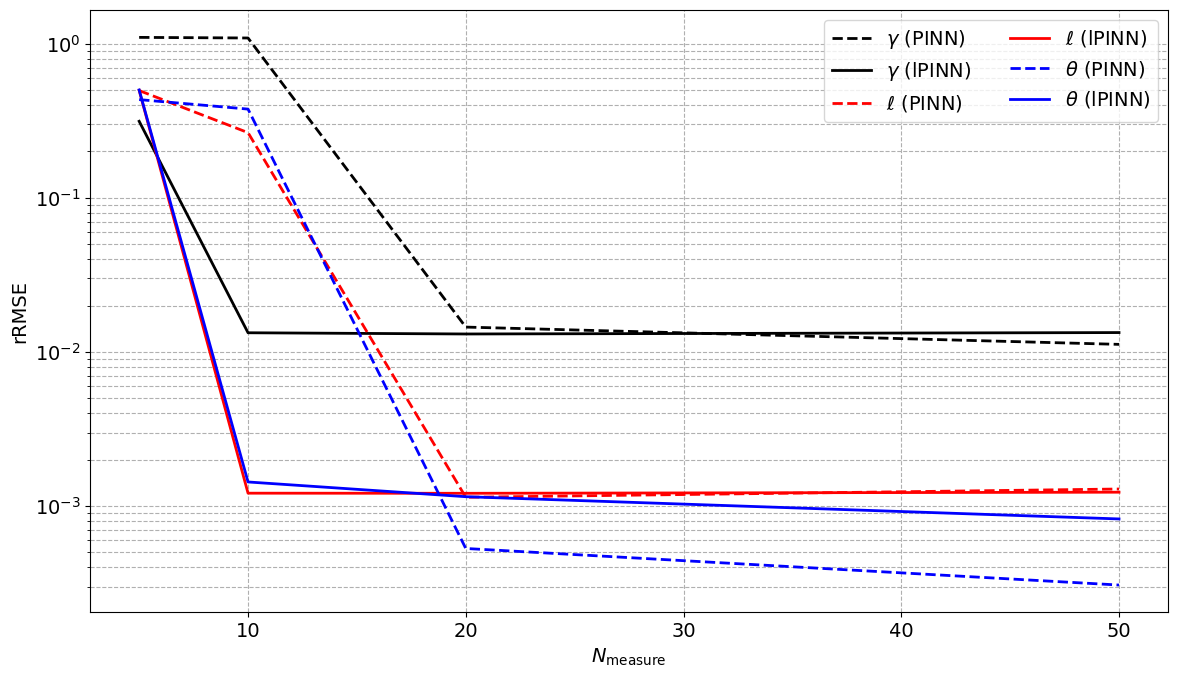}
    \caption{
    Nonlinear pendulum equation inverse problem: Relative root-mean-square error (rRMSE) of the inferred damping coefficient $\gamma$, pendulum length $\ell$, and predicted solution $\theta$ as functions of the number of measurement points $N_{\mathrm{m}}$. $N_{\mathrm{res}}=500$, $N_{\eta}=200$, $N_{\mathrm{train}}=500$.
    }
    \label{fig:Error_inv_nlpend}
\end{figure}

Figure~\ref{fig:Error_inv_nlpend} depicts the effect of the number of measurements on the accuracy of the PINN and lPINN inverse solutions.
Both PINN and lPINN errors generally decrease with increasing $N_{\mathrm{m}}$. The lPINN gives similar errors to PINN for $\theta$ and $\ell$, and estimates $\gamma$ more accurately. The lPINN errors are smaller than PINN's for $N_{\mathrm{m}} =10$ for all three unknowns. 
For $N_{\mathrm{m}}\ge 20$, lPINN gives similar errors as PINN for $\ell$ and $\gamma$, while PINN estimates $\theta$ more accurately. 

\subsection{Superresolution of lPINN}\label{sec:superresolution}

In this section, we study the interdependence between numerical errors in the training dataset and lPINN accuracy as a function of the numerical resolution of the model generating the training dataset for offline pretraining of the lPINN DNNs. The goal is to demonstrate that an lPINN pretrained on a dataset computed on a coarser mesh can produce a numerical solution that is more accurate than the test numerical solution computed on the same coarse mesh and, more importantly, is more accurate than the test numerical solution computed on a finer mesh. We refer to the latter as superresolution. 
We demonstrate the superresolution properties of lPINN for the ADE， nonlinear Burgers' and pendulum equations. 

\subsubsection{ADE equation}

We first study how lPINN accuracy depends on the accuracy of the numerical simulations that compose the training dataset. The numerical solutions are obtained using the second-order central differences in space and the adaptive fifth-order Runge–Kutta (RK45) method in time. 
We generate three training datasets, $\{ D^i \}_{i=1}^3$, each consisting of $N_{\mathrm{train}}=500$ numerical simulations, performed on the uniform mesh $(N^i_x,N^i_t)$, where  $N^i_x=N^i_t \in (30,59,240)$. We pretrained five lPINN models on each corresponding dataset. In the online stage, we use the lPINN models to obtain test solutions with $N_{\mathrm{res}}=500$ residual points. We compute the rRMSE errors with respect to the analytical solution on the mesh used to generate the corresponding dataset. For comparison, we also compute the rRMSE of the numerical test solution obtained on the corresponding mesh. 
 
\begin{figure}[H]
    \centering
    \includegraphics[width=0.75\linewidth]{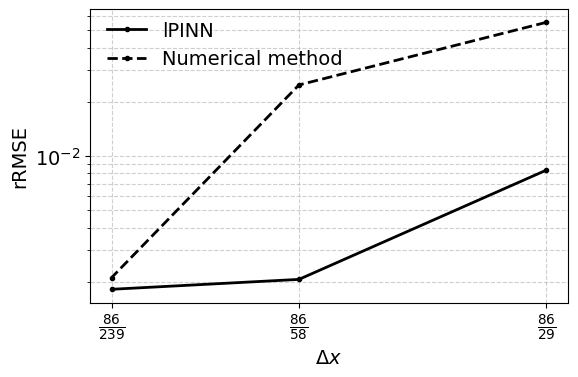}
    \caption{lPINN prediction error vs. mesh spacing $\Delta x$ for ADE equation. For each mesh, the lPINN is pretrained using the numerical-solution dataset generated on that mesh and is then evaluated on the same mesh. Smaller $\Delta x$ corresponds to a finer mesh. $N_{\mathrm{res}}=500$, $N_\eta=50$, $N_{\mathrm{train}}=500$.}
    \label{fig:ADE_lPINN_vs_dx}
\end{figure}

Figure~\ref{fig:ADE_lPINN_vs_dx} shows that lPINN and numerical solution errors decrease with decreasing $\Delta x = L/(N^i_x-1)$. Two important observations are that lPINN errors are smaller than the numerical solution error for the same $\Delta x$. More importantly, the lPINN error obtained on the grid $\frac{L}{29}$ is smaller than the error of the numerical solution obtained on the finer mesh $\frac{L}{58}$. Similarly, the lPINN error corresponding to grid size $\frac{L}{58}$ is smaller than the numerical solution error corresponding to the finer mesh $\frac{L}{239}$. We also see the limitation of the lPINN superresolution in the considered setting (the training data set size, the number of basis functions, and the number of residual points): the lPINN error pretrained on the grid $\frac{L}{29}$ is larger than the error of the numerical solution obtained on the grid $\frac{L}{239}$. 

\subsubsection{Burgers' equation}

We first study how lPINN accuracy depends on the accuracy of the numerical simulations that compose the training dataset. We obtain the numerical solutions using second-order central differences in space and an adaptive fifth-order Runge–Kutta (RK45) method in time. 
We generate three training datasets, $\{ D^i \}_{i=1}^3$, each consisting of $N_{\mathrm{train}}=500$ numerical simulations, performed on the uniform mesh $(N^i_x,N^i_t)$, where  $N^i_x=N^i_t \in (15,30,59)$. We pretrained five lPINN models on each corresponding dataset. In the online stage, we use the lPINN models to obtain test solutions with $N_{\mathrm{res}}=100$ residual points. We compute the rRMSE errors with respect to the analytical solution on the mesh used to generate the corresponding dataset. For comparison, we also compute the rRMSE of the numerical test solution obtained on the corresponding mesh. 
 
\begin{figure}[H]
    \centering
    \includegraphics[width=0.75\linewidth]{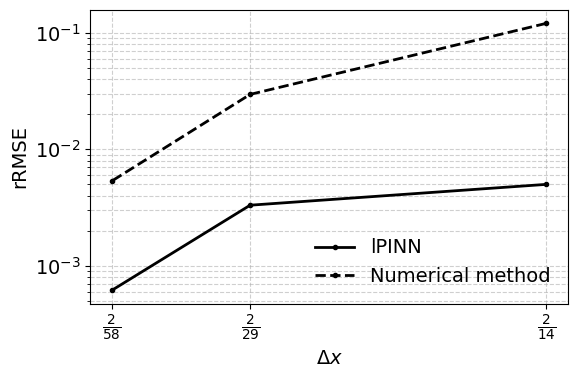}
    \caption{lPINN prediction error versus mesh spacing $\Delta x$ for Burgers' equation. For each mesh, the lPINN is pretrained using the numerical-solution dataset generated on that mesh and is then evaluated on the same mesh. Smaller $\Delta x$ corresponds to a finer mesh. $N_{\mathrm{res}}=100$, $N_\eta=20$, $N_{\mathrm{train}}=500$.}
    \label{fig:burgers_lPINN_vs_dx}
\end{figure}

Figure~\ref{fig:burgers_lPINN_vs_dx} shows that lPINN and numerical solution errors decrease with decreasing $\Delta x = L/(N^i_x-1)$. Two important observations are that lPINN errors are smaller than the numerical solution error for the same $\Delta x$. More importantly, the lPINN error obtained on the grid $\frac{L}{14}$ is smaller than the error of the numerical solution obtained on the finer mesh $\frac{L}{29}$. Similarly, the lPINN error corresponding to grid size $\frac{L}{29}$ is smaller than the numerical solution error corresponding to the finer mesh $\frac{L}{58}$. We also see the limitation of the lPINN superresolution in the considered setting (the training data set size, the number of basis functions, and the number of residual points): the lPINN error pretrained on the grid $\frac{L}{14}$ is larger than the error of the numerical solution obtained on the grid $\frac{L}{58}$.      

\subsubsection{Nonlinear Pendulum equation}
We next investigate whether the same conclusion holds for the nonlinear pendulum equation, for which the lPINN is pretrained using the derivative-matching approach. We obtain the numerical solutions using the fourth-order Runge--Kutta (RK4) method. We generate three training datasets, $\{ D^i \}_{i=1}^3$ each consisting of $N_{\mathrm{train}}=500$ numerical simulations performed on a uniform temporal grid with $N_t^i\in{(181,241,481)}$ over $t\in[0,30]$, corresponding to a constant timestep $\Delta t=T/(N_t^i-1)$. We pretrain lPINN models on each dataset. In the online stage, the lPINN models are used to obtain test solutions for $\gamma=0.1$ and $\ell=0.8$ using $N_{\mathrm{res}}=1000$ randomly sampled residual points. 
The lPINN rRMSEs are compared to the rRMSEs of the numerical solutions computed on the same meshes and with $N_t=361$.  All rRMSEs are computed with respect to a high-resolution numerical reference solution, obtained using the adaptive fifth-order Runge--Kutta (RK5) method with $N_t^{\mathrm{ref}}=2881$. 

We consider two sets of lPINN solutions, one set obtained with  $N_{\mathrm{train}} = 500$ and $N_{\eta} = 200$ and the other obtained with $N_{\mathrm{train}} = 4000$ and $N_{\eta} = 400$. 

\begin{figure}[H]
    \centering
    \includegraphics[width=0.75\linewidth]{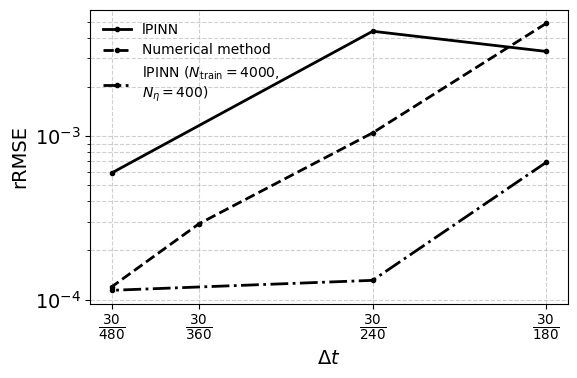}
    \caption{lPINN prediction error $e_{\mathrm{lPINN}}$ versus time-step size $\Delta t$ for the nonlinear pendulum equation. For each temporal mesh, the lPINN is pretrained using the numerical-solution dataset generated on that mesh and is then evaluated on the same mesh. Smaller $\Delta t$ corresponds to a finer temporal mesh. 
   $N_{\mathrm{res}}=1000$.}
    \label{fig:nlpend_lPINN_vs_dt}
\end{figure}

Figure~\ref{fig:nlpend_lPINN_vs_dt} plots rRMSEs of the numerical and lPINN solutions as functions of the time step size. 
With $N_{\mathrm{train}} = 500$ and $N_{\eta} = 200$, lPINN does not produce superresolution. Moreover, lPINN does not consistently outperform numerical solutions obtained on the same temporal grids. Increasing the ensemble size to $N_{\mathrm{train}} = 4000$ and the basis dimension to $N_{\eta} = 400$ substantially improves the accuracy of the basis function representation, as evidenced by the reduction in the lPINN solution rRMSEs. Under this enriched setting, the lPINN rRMSEs are smaller than the rRMSEs of the numerical solutions obtained on the same temporal meshes. Furthermore, the lPINN trained with $N_t = 181$ is more accurate than the numerical solution obtained with $N_t = 241$, and the lPINN trained with $N_t = 241$ is more accurate than the numerical solution obtained with $N_t = 361$. Thus, the lPINN with $N_{\mathrm{train}} = 4000$ and $N_{\eta} = 400$ achieves superresolution.

\section{Conclusion}

We developed the Linearized Physics-Informed Neural Network (lPINN) as a physics-corrected neural reduced-basis method for many-query forward and inverse problems governed by differential equations. The central feature of the method is that an ensemble of numerical solutions is used offline to learn a continuous, coordinate-dependent neural basis functions, while the solution of each new problem is determined online by minimizing the governing-equation residual in the resulting low-dimensional coefficient space. 
The continuous neural basis can be pretrained through derivative matching or by augmenting solution reconstruction with physics residuals. These alternatives provide different mechanisms for representing the differential quantities required during online inference. The numerical experiments show that the preferred strategy can depend on the governing equation and neural architecture. Residual-based pretraining produced the most accurate downstream solutions for the advection--diffusion and Burgers' equations, whereas derivative-matching pretraining was more effective for the oscillatory nonlinear pendulum problem. The latter example also demonstrated that the framework can incorporate Fourier feature networks when the solution family contains oscillatory components that are difficult to represent with a standard multilayer perceptron. Once the neural trial space has been constructed, lPINN restricts problem-specific optimization to the reduced coefficients and, in inverse problems, the unknown physical parameters. For linear differential operators, the online problem can be solved as regularized linear least squares. For nonlinear operators, the residual remains nonlinear in the reduced variables, but the nonlinear feature representation is not retrained. This distinction enables lPINN to extend final-layer adaptation beyond linear forward problems to nonlinear forward and inverse inference. The method was evaluated on the advection--diffusion equation, Burgers' equation, and the nonlinear pendulum equation. In the reported test cases, lPINN produced lower solution and parameter errors than the corresponding vanilla PINNs while reducing online inference times by factors ranging from approximately one order of magnitude to more than three orders of magnitude. Its advantages were generally most pronounced when the numbers of residual or measurement points were limited. The cross-resolution experiments further showed that a continuous basis learned from coarse-mesh solution data could be evaluated on finer meshes without retraining and with nearly unchanged accuracy. These results demonstrate cross-resolution transfer of the learned representation, although the attainable accuracy remains limited by the quality and expressiveness of the offline trial space. 

Several limitations define the scope of the present results. First, lPINN requires a representative ensemble of numerical solutions, and its reliability outside the parameter, initial-condition, or boundary-condition ranges covered during pretraining has not been established. Second, generating the training ensemble and learning the neural basis can be computationally expensive. The method is therefore most appropriate when this offline cost can be amortized over a sufficiently large number of forward or inverse queries. Third, performance depends on the basis dimension, architecture, pretraining objective, regularization, and online enforcement of initial and boundary conditions. The observed degradation of constraint accuracy in some experiments also indicates that properties learned offline should not always be assumed to remain satisfied when only the interior differential-equation residual is minimized online. 

Future work should establish systematic criteria for selecting the basis dimension and pretraining strategy, develop adaptive enrichment procedures for parameter regions that are not adequately represented by the initial ensemble, and investigate exact or weighted enforcement of initial and boundary conditions during reduced inference. Additional priorities include quantifying out-of-distribution generalization, analyzing amortized computational cost and conditioning, and comparing the method with classical residual-minimizing reduced-order models, random-feature methods, and neural operators. Extensions to higher-dimensional systems, spatially heterogeneous coefficients, variable geometries, and field-scale forward and inverse problems will be necessary to determine the practical range of the approach. 

Overall, the results support the use of lPINN as an alternative to repeated full-network PINN optimization in many-query settings. By learning an operator-compatible continuous neural reduced basis offline and enforcing the governing physics in its coefficient space online, lPINN combines reusable data-driven representation with problem-specific physical correction.

The present comparisons are designed to isolate the computational consequences of replacing instance-specific full-network PINN optimization with reduced coefficient inference in a pretrained neural trial space. They should not be interpreted as establishing superiority over classical or nonlinear reduced-order models. A systematic ROM comparison would require controlling for the snapshot ensemble, basis dimension, continuous reconstruction of discrete modes, derivative approximation, residual sampling, constraint enforcement, and hyper-reduction. Such a comparison is an important direction for future work, particularly for quantifying when the intrinsic coordinate continuity and automatic differentiability of the lPINN basis justify its additional offline training cost.

\section{Acknowledgements}
This research was supported by CUSSP (Center for Understanding Subsurface Signals and Permeability), an Energy Earthshot Research Center funded by the U.S. Department of Energy (DOE), Office of Science under FWP 81834, and the SRI (Strategic Research Initiative) Program at the UIUC’s Grainger College of Engineering. 

\section{Appendix}

\subsection{Summary of Hyperparameters, Error and Standard Deviation}

\begin{table}[H]
\centering
\caption{
Forward Problem in Basis Size Studies for lPINN. $N_{\mathrm{train}}=500$.
For ADE equation, $N_{\mathrm{IC}}=N_{\mathrm{res}}$ and
$N_{\mathrm{BC}}=2N_{\mathrm{res}}$.
For Burgers' equation, $N_{\mathrm{IC}}=N_{\mathrm{BC}}=N_{\mathrm{res}}$.
}
\label{tab:basis_size_hyperparameters}
\resizebox{\textwidth}{!}{
\begin{tabular}{
c l c c c c c c c
}
\toprule
Fig.
& Equation
& $N_{\eta}$
& $N_{\mathrm{res}}$
& $\lambda_{\mathrm{IC},0}$
& $\lambda_{\mathrm{BC}}$
& $\lambda$
& rRMSE$_u$
& $\sigma_{\mathrm{rRMSE}_u}$\\
\midrule

\multirow{6}{*}{2}
& \multirow{6}{*}{ADE}
& 10  & 100
& 0.7 & 1 & 0
& 0.002579 & 0.001858 \\
&
& 20  & 200
& 1.5 & 1.75 & $10^{-4}$
& 0.004899 & 0.002794 \\
&
& 50  & 500
& 0.3 & 0.5 & 0
& 0.001249 & 0.002268 \\
&
& 100 & 1000
& 0.6 & 1 & $10^{-4}$
& 0.001342 & 0.002638 \\
&
& 150 & 1500
& 0.6 & 0.85 & 0
& 0.001555 & 0.001287 \\
&
& 200 & 2000
& 0.7 & 1 & 0
& 0.001499 & 0.001388 \\

\midrule

\multirow{4}{*}{8}
& \multirow{4}{*}{Burgers'}
& 8 & 80
& 1 & 10 & 0
& 0.000962 & 0.000280
\\
&
& 20 & 200
& 10 & 1 & $10^{-8}$
& 0.000741 & 0.000185
\\
&
& 50 & 500
& 30 & 3 & 0
& 0.001526 & 0.000159
\\
&
& 100 & 1000
& 10 & 0 & 0
& 0.001880 & 0.000218
\\

\bottomrule
\end{tabular}
}
\end{table}

\begin{table}[H]
\centering
\caption{
Nonlinear Pendulum Forward Problem in the Basis Size Study for lPINN.
$N_{\mathrm{train}}=500$.
}
\label{tab:basis_size_forward_pendulum}
\begin{tabular}{
c c c c c c c c
}
\toprule
Fig.
& $N_{\eta}$
& $N_{\mathrm{res}}$
& $\lambda_{\mathrm{IC},0}$
& $\lambda_{\mathrm{IC},1}$
& $\lambda$
& rRMSE$_\theta$
& $\sigma_{\mathrm{rRMSE}_\theta}$\\
\midrule

\multirow{3}{*}{13}
& 100 & 1000
& 316.228 & 0.003 & 0
& 0.050856 & 0.005259
\\
&
200 & 2000
& 0.01 & 10 & 0
& 0.002318 & 0.003529
\\
&
400 & 4000
& 0 & 0 & $10^{-4}$
& 0.003344 & 0.000264
\\

\bottomrule
\end{tabular}
\end{table}

\begin{table}[H]
\centering
\caption{
ADE Inverse Problem in the Basis Size Study for lPINN. $N_{\mathrm{train}}=500$, $N_{\mathrm{m}}=40$, $N_{\mathrm{IC}}=N_{\mathrm{BC}}=N_{\mathrm{res}}$, $\lambda_{\mathrm{data}}=1$, Fig.~2.
}
\label{tab:basis_size_inverse_ade}
\resizebox{\textwidth}{!}{
\begin{tabular}{
c c c
c c
c c
c c
}
\toprule
$N_{\eta}$
& $N_{\mathrm{res}}$
& $\lambda$
& rRMSE$_u$
& $\sigma_{\mathrm{rRMSE}_u}$
& rRMSE$_V$& $\sigma_{\mathrm{rRMSE}_V}$
& rRMSE$_D$& $\sigma_{\mathrm{rRMSE}_D}$
\\
\midrule

10  & 100  & 0
& 0.007071 & 0.004013
& 0.003085 & 0.002319
& 0.194328 & 0.221529
\\

20  & 200  & 0
& 0.008056 & 0.000843
& 0.005200 & 0.001114
& 0.197752 & 0.035221
\\

50  & 500  & $10^{-6}$
& 0.003898 & 0.002181
& 0.004135 & 0.002421
& 0.107233 & 0.088383
\\

100 & 1000 & $10^{-8}$
& 0.006063 & 0.001795
& 0.004601 & 0.001695
& 0.166546 & 0.067861
\\

150 & 1500 & 0
& 0.005475 & 0.001152
& 0.003259 & 0.001233
& 0.134873 & 0.046461
\\

200 & 2000 & $10^{-6}$
& 0.006883 & 0.001085
& 0.005783 & 0.000727
& 0.195849 & 0.046107
\\

\bottomrule
\end{tabular}
}
\end{table}

\begin{table}[H]
\centering
\caption{
Burgers Equation Inverse Problem in the Basis Size Study for lPINN. $N_{\mathrm{train}}=500$, $N_{\mathrm{m}}=50$, $N_{\mathrm{IC}}=N_{\mathrm{BC}}=N_{\mathrm{res}}$, $\lambda_{\mathrm{data}}=1$, Fig.~8.
}
\label{tab:basis_size_inverse_burgers}
\begin{tabular}{
c c c
c c
c c
}
\toprule
$N_{\eta}$
& $N_{\mathrm{res}}$
& $\lambda$
& rRMSE$_u$
& $\sigma_{\mathrm{rRMSE}_u}$
& rRMSE$_\nu$
& $\sigma_{\mathrm{rRMSE}_\nu}$
\\
\midrule

8 & 80 & 0
& 0.005800 & 0.001245
& 0.039483 & 0.008887
\\

20 & 200 & $10^{-6}$
& 0.006070 & 0.002712
& 0.042535 & 0.019688
\\

50 & 500 & $10^{-8}$
& 0.002922 & 0.000346
& 0.014586 & 0.003803
\\

100 & 1000 & $10^{-8}$
& 0.003454 & 0.000325
& 0.018585 & 0.002970
\\

\bottomrule
\end{tabular}
\end{table}

\begin{table}[H]
\centering
\caption{
Nonlinear Pendulum Inverse Problem in the Basis Size Study for lPINN.
$N_{\mathrm{train}}=500$, $N_{\mathrm{m}}=50$,  $\lambda_{\mathrm{data}}=1$,  Fig.~13.
}
\label{tab:basis_size_inverse_pendulum}
\resizebox{\textwidth}{!}{
\begin{tabular}{
c c c
c c
c c
c c
}
\toprule
$N_{\eta}$
& $N_{\mathrm{res}}$
& $\lambda$
& rRMSE$_\theta$
& $\sigma_{\mathrm{rRMSE}_\theta}$
& rRMSE$_\gamma$& $\sigma_{\mathrm{rRMSE}_\gamma}$
& rRMSE$_\ell$& $\sigma_{\mathrm{rRMSE}_\ell}$\\
\midrule

100 & 1000 & 0
& 0.006772 & 0.000292
& 0.04747 & 0.006382
& 0.00092 & 0.000347
\\

200 & 2000 & $10^{-8}$
& 0.001398 & 0.000245
& 0.02072 & 0.001825
& 0.00083 & 0.000248
\\

400 & 4000 & $10^{-8}$
& 0.000600 & 0.000028
& 0.01841 & 0.000628
& 0.00134 & 0.000025
\\

\bottomrule
\end{tabular}
}
\end{table}

\begin{table}[H]
\centering
\caption{
Forward Problem: PINN and lPINN Error vs.\ $N_{\mathrm{res}}$.
$N_{\mathrm{train}}=500$ for lPINN.
For ADE equation, $N_{\mathrm{IC}}=N_{\mathrm{res}}$ and
$N_{\mathrm{BC}}=2N_{\mathrm{res}}$.
For Burgers' equation,
$N_{\mathrm{IC}}=N_{\mathrm{BC}}=N_{\mathrm{res}}$.
}
\label{tab:forward_hyperparameters}

\resizebox{\textwidth}{!}{
\begin{tabular}{
c l r l
c c c c
c c
}
\toprule
Fig.
& Equation
& $N_{\mathrm{res}}$
& Method
& $N_{\eta}$
& $\lambda_{\mathrm{IC},0}$
& $\lambda_{\mathrm{BC}}$
& $\lambda$
& rRMSE$_u$
& $\sigma_{\mathrm{rRMSE}_u}$
\\
\midrule

4
& \multirow{8}{*}{ADE}
& 100
& PINN & --
& 1 & 1 & $10^{-4}$
& 0.028976 & 0.026853
\\

4
&
& 100
& lPINN & 50
& 0.003 & 0.01 & $10^{-4}$
& 0.002354 & 0.005093
\\

3,4
&
& 500
& PINN & --
& 1 & 1 & $10^{-4}$
& 0.010406 & 0.004776
\\

3,4
&
& 500
& lPINN & 50
& 0.1 & 0.1 & $10^{-4}$
& 0.001451 & 0.002036
\\

4
&
& 1500
& PINN & --
& 1 & 1 & $10^{-4}$
& 0.010727 & 0.004047
\\

4
&
& 1500
& lPINN & 50
& 0.3 & 0.3 & $10^{-4}$
& 0.001420 & 0.001167
\\

4
&
& 3000
& PINN & --
& 1 & 1 & $10^{-4}$
& 0.010341 & 0.002033
\\

4
&
& 3000
& lPINN & 50
& 3.225 & 6.45 & $10^{-4}$
& 0.001447 & 0.000705
\\

\midrule


9,10
& \multirow{8}{*}{Burgers'}
& 100
& PINN & --
& 1 & 2 & $10^{-6}$
& 0.043862 & 0.099357
\\

9,10
&
& 100
& lPINN & 20
& 24 & 0.03 & $10^{-6}$
& 0.000807 & 0.000189
\\

10
&
& 300
& PINN & --
& 1 & 2 & $10^{-6}$
& 0.000902 & 0.000905
\\

10
&
& 300
& lPINN & 20
& 25 & 1.00 & $10^{-6}$
& 0.000780 & 0.000047
\\

10
&
& 800
& PINN & --
& 1 & 2 & $10^{-6}$
& 0.000211 & 0.000331
\\

10
&
& 800
& lPINN & 20
& 25 & 0.15 & $10^{-6}$
& 0.000810 & 0.000041
\\

10
&
& 1500
& PINN & --
& 1 & 2 & $10^{-6}$
& 0.000151 & 0.000099
\\

10
&
& 1500
& lPINN & 20
& 25 & 0.005 & $10^{-6}$
& 0.000835 & 0.000051
\\

\bottomrule
\end{tabular}
}
\end{table}

\begin{table}[H]
\centering
\caption{
Nonlinear Pendulum Forward Problem: PINN and lPINN Error vs.\ $N_{\mathrm{res}}$.
$N_{\mathrm{train}}=500$ for lPINN.
}
\label{tab:forward_pendulum_hyperparameters}

\resizebox{\textwidth}{!}{
\begin{tabular}{
c r l
c c c c
c c
}
\toprule
Fig.
& $N_{\mathrm{res}}$
& Method
& $N_{\eta}$
& $\lambda_{\mathrm{IC},0}$
& $\lambda_{\mathrm{IC},1}$
& $\lambda$
& rRMSE$_\theta$
& $\sigma_{\mathrm{rRMSE}_\theta}$
\\
\midrule

15
& 500
& PINN & --
& 1 & 1 & $10^{-8}$
& 0.271213 & 0.112946
\\

15
& 500
& lPINN & 200
& $10^{-4}$ & 0 & $10^{-8}$
& 0.005755 & 0.013405
\\

14,15
& 1000
& PINN & --
& 1 & 1 & $10^{-8}$
& 0.099410 & 0.159099
\\

14,15
& 1000
& lPINN & 200
& 3 & 3 & $10^{-8}$
& 0.002208 & 0.003536
\\

15
& 2000
& PINN & --
& 1 & 1 & $10^{-8}$
& 0.015899 & 0.000640
\\

15
& 2000
& lPINN & 200
& 10 & 100 & $10^{-8}$
& 0.002178 & 0.004497
\\

15
& 8000
& PINN & --
& 1 & 1 & $10^{-8}$
& 0.016420 & 0.000193
\\

15
& 8000
& lPINN & 200
& 0.1 & 300 & $10^{-8}$
& 0.002187 & 0.002994
\\

\bottomrule
\end{tabular}
}
\end{table}

\begin{table}[H]
\centering
\caption{
ADE Inverse Problem: PINN and lPINN Error vs.\ $N_{\mathrm{m}}$. $N_{\mathrm{train}}=500$, $N_{\eta}=50$, $N_{\mathrm{res}}=500$, $N_{\mathrm{IC}}=500$, $N_{\mathrm{BC}}=500$, $\lambda_{\mathrm{data}}=1$. 
}
\label{tab:inverse_ade}
\resizebox{\textwidth}{!}{
\begin{tabular}{
c c l
c c llc
c c
c c
c c
}
\toprule
Fig.
& $N_{\mathrm{m}}$
& Method
& $N_{\mathrm{res}}$
& $N_{\eta}$
&  $\lambda_{\mathrm{IC},0}$& $\lambda_{\mathrm{BC}}$&$\lambda$
& rRMSE$_u$
& $\sigma_{\mathrm{rRMSE}_u}$
& rRMSE$_V$& $\sigma_{\mathrm{rRMSE}_V}$
& rRMSE$_D$& $\sigma_{\mathrm{rRMSE}_D}$
\\
\midrule

\multirow{2}{*}{6}
& 5
& PINN & 500 & 50
&  1& 1&$10^{-4}$
& 0.253604 & 0.235495
& 0.486947 & 0.513296
& 0.996616 & 0.001698
\\
&
5
& lPINN & 500 & 50
&  0& 0&0
& 0.010980 & 0.004001
& 0.006302 & 0.001981
& 0.366958 & 0.158113
\\

\multirow{2}{*}{5,6}
& 40
& PINN & 500 & 50
&  1& 1&0
& 0.050878 & 0.030874
& 0.028985 & 0.019506
& 0.851279 & 0.279618
\\
&
40
& lPINN & 500 & 50
&  0& 0&0
& 0.003886 & 0.002205
& 0.004118 & 0.002396
& 0.107385 & 0.089965
\\

\multirow{2}{*}{6}
& 80
& PINN & 500 & 50
&  1& 1&0
& 0.018262 & 0.040573
& 0.024940 & 0.029324
& 0.491649 & 0.416658
\\
&
80
& lPINN & 500 & 50
&  0& 0&0
& 0.002356 & 0.000988
& 0.002744 & 0.001640
& 0.052280 & 0.056886
\\

\multirow{2}{*}{6}
& 160
& PINN & 500 & 50
&  1& 1&$10^{-6}$
& 0.002576 & 0.004371
& 0.001215 & 0.007101
& 0.172043 & 0.281317
\\
&
160
& lPINN & 500 & 50
&  0& 0&0
& 0.001833 & 0.000939
& 0.001649 & 0.001257
& 0.051348 & 0.056867
\\

\bottomrule
\end{tabular}
}
\end{table}

\begin{table}[H]
\centering
\caption{
Burgers Equation Inverse Problem: PINN and lPINN Error vs.\
 $N_{\mathrm{m}}$. $N_{\mathrm{train}}=500$, $N_{\eta}=50$, $N_{\mathrm{res}}=100$, $N_{\mathrm{IC}}=100$, $N_{\mathrm{BC}}=100$. 
}
\label{tab:inverse_burgers}
\begin{tabular}{
c c l
c llc
c c
c c
}
\toprule
Fig.
& $N_{\mathrm{m}}$
& Method
& $\lambda_{\mathrm{data}}$
 & $\lambda_{\mathrm{IC},0}$&$\lambda_{\mathrm{BC}}$& $\lambda$
& rRMSE$_u$
& $\sigma_{\mathrm{rRMSE}_u}$
& rRMSE$(\bar{\nu})$& $\sigma_{\mathrm{rRMSE}_\nu}$
\\
\midrule

\multirow{2}{*}{12}
& 1
& PINN & 1 & 1&1& $10^{-8}$
& 0.089693 & 0.040345
& 0.43077 & 0.193140
\\
&
1
& lPINN & 1 & 0&0& $10^{-8}$
& 0.002252 & 0.001035
& 0.00786 & 0.006911
\\

\multirow{2}{*}{11,12}
& 5
& PINN & 1 & 1&1& $10^{-8}$
& 0.006089 & 0.016743
& 0.01190 & 0.027147
\\
&
5
& lPINN & 1 & 0&0& $10^{-6}$
& 0.001568 & 0.000236
& 0.00318 & 0.002070
\\

\multirow{2}{*}{12}
& 20
& PINN & 1 & 1&1& 0
& 0.001236 & 0.001589
& 0.00214 & 0.011031
\\
&
20
& lPINN & 5  & 0&0& $10^{-8}$
& 0.001503 & 0.000098
& 0.00196 & 0.001913
\\

\multirow{2}{*}{12}
& 50
& PINN & 1 & 1&1& $10^{-6}$
& 0.001059 & 0.000894
& 0.00173 & 0.004184
\\
&
50
& lPINN & 1000  & 0&0& $10^{-8}$
& 0.001409 & 0.000052
& 0.001456 & 0.001385
\\

\bottomrule
\end{tabular}
\end{table}

\begin{table}[H]
\centering
\caption{
Nonlinear Pendulum Inverse Problem: PINN and lPINN Error vs.\
$N_{\mathrm{m}}$. $N_{\mathrm{train}}=500$, $N_{\eta}=200$, $N_{\mathrm{res}}=500$. 
}
\label{tab:inverse_pendulum}
\resizebox{\textwidth}{!}{
\begin{tabular}{
c c l
c llc
c c
c c
c c
}
\toprule
Fig.
& $N_{\mathrm{m}}$
& Method
& $\lambda_{\mathrm{data}}$
 & $\lambda_{\mathrm{IC},0}$&$\lambda_{\mathrm{IC},1}$& $\lambda$
& rRMSE$_\theta$
& $\sigma_{\mathrm{rRMSE}_\theta}$
& rRMSE$_\gamma$& $\sigma_{\mathrm{rRMSE}_\gamma}$
& rRMSE$_\ell$& $\sigma_{\mathrm{rRMSE}_\ell}$
\\
\midrule

\multirow{2}{*}{17}
& 5
& PINN & 1 & 10&10& $10^{-4}$
& 0.434149 & 0.048105
& 1.10 & 0.337883
& 0.497 & 0.107659
\\
&
5
& lPINN & 1 & 0&0& $10^{-4}$
& 0.501447 & 0.000623
& 0.314 & 0.018098
& 0.499 & 0.000960
\\

\multirow{2}{*}{16,17}
& 10
& PINN & 1 & 10&10& 0
& 0.376713 & 0.135898
& 1.09 & 0.507195
& 0.265 & 0.161166
\\
&
10
& lPINN & 5  & 0&0& $10^{-8}$
& 0.001433 & 0.000873
& 0.01332 & 0.006621
& 0.001212 & 0.000433
\\

\multirow{2}{*}{17}
& 20
& PINN & 1 & 10&10& $10^{-4}$
& 0.000530 & 0.000107
& 0.0145 & 0.003320
& 0.00114 & 0.000068
\\
&
20
& lPINN & 5  & 0&0& $10^{-8}$
& 0.001150 & 0.000732
& 0.01308 & 0.009210
& 0.00121 & 0.000390
\\

\multirow{2}{*}{17}
& 50
& PINN & 1 & 10&10& $10^{-6}$
& 0.000307 & 0.000244
& 0.0112 & 0.005170
& 0.00129 & 0.000288
\\
&
50
& lPINN & 1000  & 0&0& $10^{-8}$
& 0.000825 & 0.000112
& 0.01336 & 0.006976
& 0.00123 & 0.000387
\\

\bottomrule
\end{tabular}
}
\end{table}

\begin{table}[H]
\centering
\caption{
ADE lPINN and dPICKLE Comparison.
$N_{\mathrm{train}}=500$, $N_{\eta}=50$,
$N_{\mathrm{res}}=784$, $N_{\mathrm{IC}}=30$, and
$N_{\mathrm{BC}}=59$.
}
\label{tab:dpickle_hyperparameters}
\begin{tabular}{
l c c c c
}
\toprule
Method
& $\lambda_{\mathrm{IC},0}$
& $\lambda_{\mathrm{BC}}$
& $\lambda$
& rRMSE$_u$
\\
\midrule

lPINN
& 0.5
& 750
& $10^{-4}$
& 0.0013
\\

dPICKLE
& --
& --
& --
& 0.0470
\\

\bottomrule
\end{tabular}
\end{table}

\begin{table}[H]
\centering
\caption{
Superresolution studies.
For the ADE equation,
$N_{\mathrm{IC}}=N_{\mathrm{res}}$ and
$N_{\mathrm{BC}}=2N_{\mathrm{res}}$.
For Burgers' equation,
$N_{\mathrm{IC}}=N_{\mathrm{BC}}=N_{\mathrm{res}}$. $N_{\mathrm{train}}=500$ for lPINN, 
rRMSE$_{\mathrm{num}}$ is the error of the numerical solution.
}
\label{tab:mesh_hyperparameters}
\resizebox{\textwidth}{!}{
\begin{tabular}{
c l c c c c c c c c c
}
\toprule
Fig.
& Equation
& Mesh ($N_t \times N_x$)
& $N_{\eta}$
& $N_{\mathrm{res}}$
& $\lambda_{\mathrm{IC},0}$
& $\lambda_{\mathrm{BC}}$
& $\lambda$
& rRMSE$_{\mathrm{num}}$
& rRMSE$_u$
& $\sigma_{\mathrm{rRMSE}_u}$
\\
\midrule


\multirow{3}{*}{18}
& \multirow{3}{*}{ADE}
& $30\times30$
& 50 & 500
& 7 & 6 & 0
& 0.055363
& 0.008344 & 0.002795
\\
&
& $59\times59$
& 50 & 500
& 3.5 & 5 & $10^{-8}$
& 0.024839
& 0.002061 & 0.002378
\\
&
& $240\times240$
& 50 & 500
& 0.48 & 0.96 & 0
& 0.002109
& 0.001816 & 0.001932
\\

\midrule


\multirow{3}{*}{19}
& \multirow{3}{*}{Burgers'}
& $15\times15$
& 20 & 100
& 30 & 70 & 0
& 0.120381
& 0.005004 & 0.000806
\\
&
& $30\times30$
& 20 & 100
& 8 & 0.1 & $10^{-4}$
& 0.029587
& 0.003316 & 0.000202
\\
&
& $59\times59$
& 20 & 100
& 0.3 & 10 & $10^{-8}$
& 0.005350
& 0.000618 & 0.000088
\\

\bottomrule
\end{tabular}
}
\end{table}

\begin{table}[H]
\centering
\caption{
Nonlinear Pendulum Superresolution Study. $N_{\mathrm{res}}=1000$ for lPINN. 
rRMSE$_{\mathrm{num}}$ is the error of the numerical solution.
}
\label{tab:mesh_hyperparameters_pendulum}
\resizebox{\textwidth}{!}{
\begin{tabular}{
c c c c c c c c c c
}
\toprule
Fig.
& $N_t$
& $N_{\mathrm{train}}$
& $N_{\eta}$
& $\lambda_{\mathrm{IC},0}$
& $\lambda_{\mathrm{IC},1}$
& $\lambda$
& rRMSE$_{\mathrm{num}}$
& rRMSE$_\theta$
& $\sigma_{\mathrm{rRMSE}_\theta}$
\\
\midrule

\multirow{6}{*}{20}
& 181
& 500 & 200 & 10 & 10 & $10^{-8}$
& 0.004882
& 0.003288 & 0.002839
\\

&
241
& 500 & 200 & 10 & 1 & 0
& 0.001046
& 0.004369 & 0.003157
\\

&
481
& 500 & 200 & 100 & 100 & $10^{-6}$
& 0.000120
& 0.000597 & 0.001220
\\

&
181
& 4000 & 400 & 10 & 10 & 0
& 0.004882
& 0.000691 & 0.001927
\\

&
241
& 4000 & 400 & 0.1 & 1 & $10^{-6}$
& 0.001046
& 0.000131 & 0.000760
\\

&
481
& 4000 & 400 & 0.1 & 0.1 & 0
& 0.000120
& 0.000114 & 0.000362
\\

\bottomrule
\end{tabular}
}
\end{table}

\begin{table}[H]
    \centering
    \caption{
    Offline pretraining. For ADE and Burgers' equation, $\lambda_f=\lambda_{\mathrm{IC},0}=\lambda_{\mathrm{BC}}=1$ in \ref{eq:joint_train_physics}, and for nonlinear pendulum equation, $\lambda_1=\lambda_4=1$ in  \ref{eq:nlpend_mean_pretrain} and  \ref{eq:nlpend_fluctuation_pretrain}. 
    }
    \label{tab:offline_pretraining_weights}

    \renewcommand{\arraystretch}{1.2}
    \setlength{\tabcolsep}{6pt}

    \begin{tabular}{
        l
        c
    }
        \toprule
        Equation
        & $\lambda$\\
        \midrule

        ADE
        & $10^{-6}$
        \\

        Burgers
        & $10^{-6}$
        \\

        Nonlinear pendulum mean network& 0
        \\

        Nonlinear pendulum basis& $10^{-6}$
        \\

        \bottomrule
    \end{tabular}
\end{table}

\bibliographystyle{unsrt}   
\bibliography{references}

\end{document}